\documentclass[twoside,11pt]{article}
\usepackage{jair, theapa, rawfonts}
\usepackage{xcolor}
\usepackage{graphicx}
\usepackage{mathrsfs}
\usepackage{multirow}
\usepackage{diagbox}
\usepackage{placeins}

\usepackage[ruled,vlined, linesnumbered]{algorithm2e}

\jairheading{77}{2023}{1489-1538}{01/2023}{08/2023}
\ShortHeadings{Improved Peel-and-Bound}
{Rudich, Cappart, \& Rousseau}
\firstpageno{1}

\newcommand\edit[1]{\textcolor{black}{#1}}

\begin{document}

\title{Improved Peel-and-Bound: Methods for Generating Dual Bounds with Multivalued Decision Diagrams}

\author{\name Isaac Rudich \email isaac.rudich@polymtl.ca \\
        \name Quentin Cappart \email quentin.cappart@polymtl.ca \\
        \name Louis-Martin Rousseau \email louis-martin.rousseau@polymtl.ca \\
        \addr 2500 Chem. de Polytechnique, \\Montréal, QC H3T 1J4, Canada
       }

\bibliographystyle{theapa}

\definecolor{safeA2}{RGB}{212,17,89}

\maketitle

\begin{abstract}
Decision diagrams are an increasingly important tool in cutting-edge solvers for discrete optimization. However, the field of decision diagrams is relatively new, and is still incorporating the library of techniques that conventional solvers have had decades to build. We drew inspiration from the \emph{warm-start} technique used in conventional solvers to address one of the major challenges faced by decision diagram based methods. Decision diagrams become more useful the wider they are allowed to be, but also become more costly to generate, especially with large numbers of variables. In the original version of this paper, we presented a method of \emph{peeling} off a sub-graph of previously constructed diagrams and using it as the initial diagram for subsequent iterations that we call \emph{peel-and-bound}. We tested the method on the \emph{sequence ordering problem}, and our results indicate that our \emph{peel-and-bound} scheme generates stronger bounds than a branch-and-bound scheme using the same propagators, and at significantly less computational cost. In this extended version of the paper, we also propose new methods for using relaxed decision diagrams to improve the solutions found using restricted decision diagrams, discuss the heuristic decisions involved with the parallelization of \emph{peel-and-bound}, and discuss how peel-and-bound can be hyper-optimized for sequencing problems. Furthermore, we test the new methods on the \emph{sequence ordering problem} and the \emph{traveling salesman problem with time-windows} (TSPTW), and include an updated and generalized implementation of the algorithm capable of handling any discrete optimization problem. The new results show that \emph{peel-and-bound} outperforms ddo (a decision diagram based \emph{branch-and-bound} solver) on the TSPTW. \edit{We also close 15 open benchmark instances of the TSPTW.} 
\end{abstract}

\section{Introduction}
\label{sec:introduction}
Multivalued decision diagrams (MDDs) are a useful graphical tool for compactly storing the solution space of discrete optimization problems. In the last few years, a staggering number of new applications for MDDs have been proposed~\shortcite{ddsurvey}, such as representing global constraints~\shortcite{global1,global2,global3}, handling stochastic variables~\shortcite{stochastic1,stochastic2}, and performing post-optimality analysis~\shortcite{postopt}. MDDs are particularly useful for generating strong dual 
bounds~\shortcite{cappart2019improving,dual2,dual1,dual3,cappart2022improving,gentzel2020haddock}, especially on optimization problems where linear relaxations perform poorly. There is a subset of MDD research that uses a highly paralellizable \emph{branch-and-bound} algorithm based on decision diagrams~\shortcite{BnB,ddo,bnb3,parjadis2021improving} to maximize the utility of MDD based relaxations. This paper furthers the work on the decision diagram based branch-and-bound by introducing a method, referred to as \textit{peel-and-bound}, of reusing the graphs generated at each iteration of the algorithm. This paper is an extended version of \textit{Peel-and-Bound: Generating Stronger Relaxed Bounds with Multivalued Decision Diagrams} \cite{Rudich2022PeelAndBoundGS}, the contributions of which were as follows: (1) we presented the \emph{peel-and-bound} algorithm, (2) we identified several heuristic decisions that can be used to adjust peel-and-bound, and discussed their implications, (3) we showed that peel-and-bound outperforms branch-and-bound on the \emph{sequence ordering problem} (SOP), and (4) we provided insight into how the algorithm can be applied to other problems. In this extended version of the paper we (1) provide a new implementation of the solver\footnote{https://github.com/IsaacRudich/ImprovedPnB} that is generic, and also includes an optimized framework for solving sequencing problems, (2) we propose a new heuristic for node selection, and two new search procedures that leverage the structure of peel-and-bound, (3) we test our node selection heuristics on the SOP, and finally (4) we test the new implementation on the \textit{traveling salesman problem with time windows} (TSPTW) and \textit{makespan} problem, show that it outperforms ddo (a decision diagram based branch and bound solver) \shortcite{ddo}, and (5) close several benchmark instances.

This paper is structured as follows. Section \ref{sec:background} provides the necessary technical background information and notation, as well as implementation details for the decision diagram relaxations used in our experiments. In Section \ref{sec:pnb}, we introduce the core contribution, namely the peel-and-bound procedure. The algorithm is presented, and its limitations are discussed (this section has been significantly expanded). The original computational experiments are proposed and discussed in Section \ref{sec:results}. The new computational experiments and improved implementation of the algorithm are discussed in Section \ref{sec:experiments2}.

\section{Technical Background}
\label{sec:background}

The idea of using multivalued decision diagrams (MDDs) for optimization problems was introduced by Andersen et al. (2007), and generalized by Hadzic et al. (2008) and Hoda et al. (2010). \edit{The use of decision diagrams to generate relaxed bounds was introduced by Bergman et al. (2014a).}
Following those papers, Bergman et al. (2014\edit{b}, 2016b) demonstrated the potential for a decision diagram based branch-and-bound solver to be effective, and provided an efficient parallelization scheme. Gillard et al. (2021) further improved the decision diagram based branch-and-bound solver by adding pruning techniques that can be used while the decision diagrams are being constructed, as well as to remove nodes from the branch-and-bound queue. 

This paper presents a new peel-and-bound scheme for combining restricted and relaxed decision diagrams to find exact solutions. This section provides the required technical background on  how decision diagrams can be used to model sequencing problems, and how to construct restricted/relaxed diagrams. It also introduces the notations used in this paper, and details the existing algorithms considered in our experiments. 
 
\subsection{Decision Diagrams (DDs)}
\label{sec:background:dd}
Let $\mathcal{P}$ be an instance of a discrete minimization problem with $n$ variables $\{x_1,...,x_n\}$, let $Sol(\mathcal{P})$ be the set of feasible solutions to $\mathcal{P}$, let $z^*(\mathcal{P})$ be an optimal solution to $\mathcal{P}$, and let $D(x_i)$ be the domain of variable $x_i, i \in \{1, ..., n\}$. Let $\mathcal{M}$ be a multivalued decision diagram that contains potential solutions to $\mathcal{P}$. $\mathcal{M}$ is a directed acyclic graph divided into $n+1$ layers; let $\ell_u$ be the index of the layer containing node $u, u \in \mathcal{M}$, and let $L_i$ be the set containing the nodes on layer $i$. Layer $1$ contains only a root node $r$ (with no \emph{in} arcs), and layer $n+1$ contains just a terminal node $t$ (with no \emph{out} arcs). Each arc $a_{uv} \in \mathcal{M}$ goes from a node $u$ on layer $\ell_u \in \{1,...,n\}$ to a node $v$ on layer $\ell_{u+1}$ ($\ell_{u+1} = \ell_v$). Each arc $a_{uv}$ has a label representing the assignment of variable $x_{\ell_u}$ to $l \in D(x_{\ell_u})$. An arc $a_{uv}$ with label $l$ ($a_{uv} \rightarrow l$) also has a value $v(a_{uv})$ equal to the value of being at node $u$ and assigning $x_{\ell_u}$ to $l$ ($x_{\ell_u} = l$). For simplicity, we sometime refer to $v(a_{uv})$ as $v(a)$. Thus, each path from $r$ to $t$ represents the assignment of the $n$ variables to values, and a potential solution to $\mathcal{P}$.

Let $Sol(\mathcal{M})$ be the set of all paths in $\mathcal{M}$ from $r$ to $t$, and let $T^*(u)$ be the value of the shortest path from $r$ to a node $u$. If $Sol(\mathcal{M}) = Sol(\mathcal{P})$, then $\mathcal{M}$ perfectly represents the solution space of $\mathcal{P}$, and we call $\mathcal{M}$ \emph{exact}. If $\mathcal{M}$ is exact, then the value of the shortest path through the diagram is $z^*(\mathcal{P})$ (an optimal solution to $\mathcal{P}$). Let the shortest path through $\mathcal{M}$ be $z^*(\mathcal{M})$. If $Sol(\mathcal{M}) \subseteq Sol(\mathcal{P})$, then $\mathcal{M}$ represents only feasible solutions to $\mathcal{P}$, but does not necessarily represent all feasible solutions to $\mathcal{P}$. In this case, we call $\mathcal{M}$ \emph{restricted}, and use the notation $\overline{\mathcal{M}}$ to mean that $\mathcal{M}$ is restricted (a restricted diagram yields an upper bound when minimizing). The shortest path through $\overline{\mathcal{M}}$ is not guaranteed to be optimal, but it is guaranteed to be feasible. If $Sol(\mathcal{P}) \subseteq Sol(\mathcal{M})$, then $\mathcal{M}$ represents all of the feasible solutions to $\mathcal{P}$, but potentially represents infeasible solutions as well. In this case, we call $\mathcal{M}$ \emph{relaxed}, and use the notation $\underline{\mathcal{M}}$ to mean that $\mathcal{M}$ is relaxed. The shortest path through $\underline{\mathcal{M}}$ is guaranteed to be at least as good as $z^*(\mathcal{P})$, but is not guaranteed to be feasible.

 Constructing an exact decision diagram for $\mathcal{P}$ is often intractable for large values of $n$. Observe that having an exact decision diagram means that the solution to $\mathcal{P}$ can be read in polynomial time by recursively calculating the shortest path through $\mathcal{M}$, so creating an exact decision diagram for NP-hard problems, such as for the \emph{travelling salesperson problem} (TSP), is NP-hard as well~\shortcite{DDforO}. The focus of most research that uses decision diagrams for optimization is on the construction of $\overline{\mathcal{M}}$ and/or $\underline{\mathcal{M}}$. Let $w = w(\mathcal{M})$ be the width of the largest layer of $\mathcal{M}$. The creation of an exact decision diagram potentially leads to $w$ being an exponential function of $n$, but when creating $\overline{\mathcal{M}}$ and/or $\underline{\mathcal{M}}$, $w$ can be constrained to be any natural number, limiting the number of operations construction will take. Let $w_m$ be the largest width allowed during construction. As $w_m$ approaches the width necessary to create an exact decision diagram, $z^*(\overline{\mathcal{M}})$ and $z^*(\underline{\mathcal{M}})$ approach $z^*(\mathcal{P})$, but the number of operations necessary to construct the diagram also increases.
 
\begin{figure}[h!]
    \centering
    \begin{tabular}{|c||*{4}{c|}}\hline
        \diagbox{From}{To}&A&B&C&D\\\hline
        A & - & $8$ & $5$ & $0$ \\\hline
        B & \textcolor{safeA2}{X} & - & 5 & 8 \\\hline
        C & \textcolor{safeA2}{X} & $5$ & - & 5 \\\hline
        D & \textcolor{safeA2}{X} & \textcolor{safeA2}{X} & 1 & - \\\hline
    \end{tabular}
    \caption{Example of a SOP instance with transition costs $c_{row,col}$ (\textcolor{safeA2}{X} \edit{are global precedence constraints. A must come before everything else, and B must come before D}).}
    \label{fig:basicSOP}
\end{figure}

\begin{figure}[h!]
    \centering
    \begingroup
        \setlength{\tabcolsep}{20pt}
        \begin{tabular}{c c c}
            \includegraphics[scale=.55]{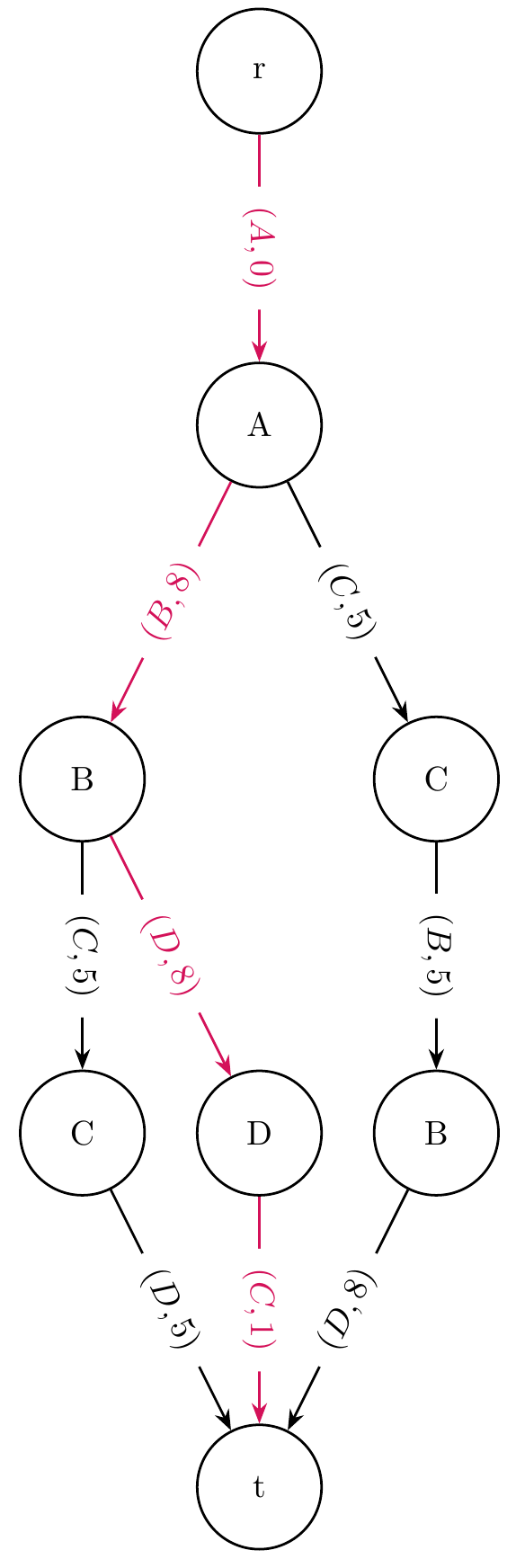} & \includegraphics[scale=.55]{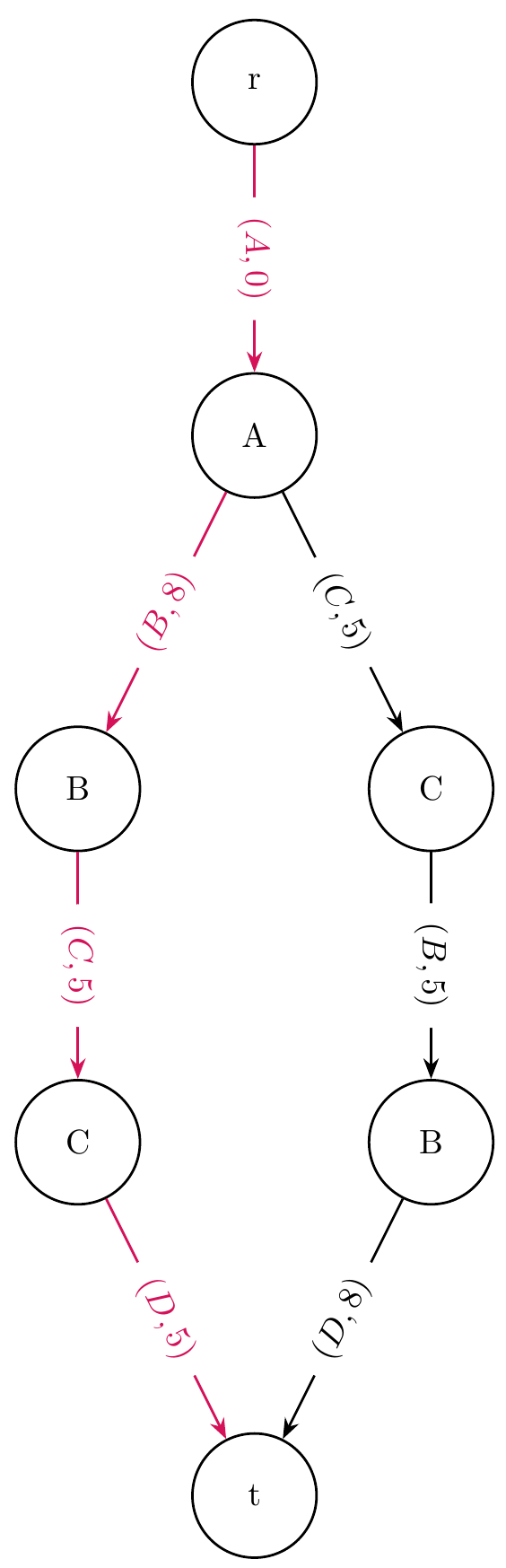} & \includegraphics[scale=.55]{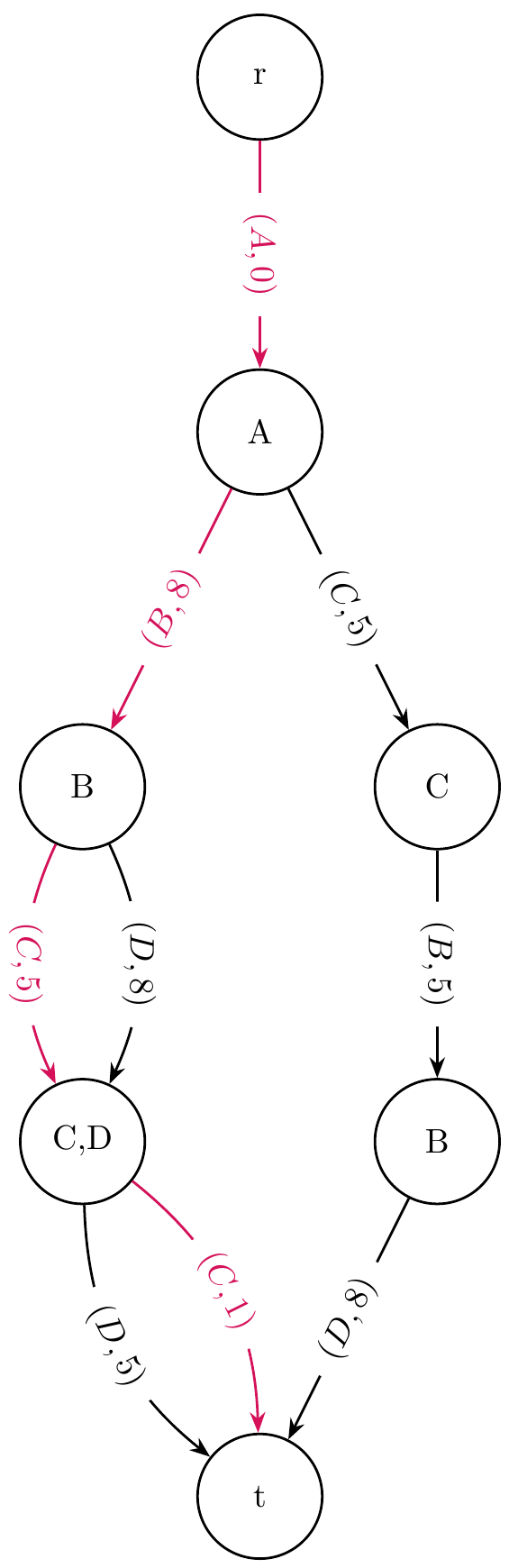}\\
            Exact Diagram & Restricted Diagram & Relaxed Diagram\\
            $z^* = {[A,B, D,C]}$ & $z^* = {[A,B, C,D]}$ & $z^* = {[A,B, C,C]}$\\
            $T^*=17$ & $T^*=18$ & $T^*=14$
        \end{tabular}
    \endgroup
    \caption{MDD Representation for the SOP instance presented in Figure \ref{fig:basicSOP}. Each arc $a$ has the format: $(l, v(a))$. The red path in each diagram indicates the shortest path from $r$ to $t$.}
    \label{fig:basicDDs}
\end{figure}

Figure \ref{fig:basicSOP} gives an instance of \emph{sequence ordering problem} (SOP), and Figure \ref{fig:basicDDs} contains simple examples of exact, restricted, and relaxed decision diagrams for that instance where $w_m=2$ for $\overline{\mathcal{M}}$ and $\underline{\mathcal{M}}$. 

The SOP requires finding the minimum-cost sequence of $n$ elements that includes each element exactly once, subject to transition costs $c_{ij}$ of following $x_i$ with $x_j$, and subject to precedence constraints requiring that certain elements precede others in the sequence. In other words, the SOP is an asymmetric TSP with precedence constraints. The label of each node matches the union of the labels of the incoming arcs. Each arc $a_{uv}$ is labeled in the format $(l, val(a_{uv}))$, representing the assignment of $x_{\ell_u}$ to $l$, and $val(a)$ represents the cost of the shortest path from the label of $u$ to $l$. In other words, an arc with label $l$ leaving layer $i$, represents the assignment of $l$ to the $i^{th}$ position of the sequence. The red path in each diagram indicates the shortest path through the diagram, and $T^*$ indicates the cost of the shortest path through the diagram.

\subsection{Restricted Decision Diagrams}
\label{sec:background:restricted}
Constructing $\overline{\mathcal{M}}$ for a given width $w_m$ is a straightforward process that can be thought of as a generalized greedy algorithm. Beginning with the root node $r$, an arc is generated for every element in the domain of $r$, and a node is generated at the end of each arc in the second layer. The process is repeated for each layer, except layer $n$ where all outgoing arcs point to the terminal, unless $w(\overline{\mathcal{M}})$ exceeds $w_m$. Then the least promising node is removed from the offending layer until $w(\overline{\mathcal{M}})$ is equal to $w_m$. The definition of \emph{least promising} is a heuristic decision. For the purposes of this paper, the least promising node is the node $u$ such that the shortest path from $r$ to $u$ is longer than the shortest path from $r$ to any other node $v\neq u$ in layer $\ell_u$. \edit{In other words, nodes with the worst associated partial solutions are chosen for removal.}

It is of note that another method of reducing the width of $\overline{\mathcal{M}}$ is merging equivalent nodes. In the SOP, two nodes can be considered equivalent if they have the same state (last element in the sequence), and all incoming paths have visited the same set of elements. For example, a node with exactly one incoming path $[A,B,C]$ could be merged with a node in the same layer with exactly one incoming path $[B,A,C]$. In many MDD applications this is a valuable insight, and it helps motivate the algorithm for constructing $\underline{\mathcal{M}}$. However, for the SOP, we observed that the work of finding equivalent nodes in $\overline{\mathcal{M}}$ often outweighed the benefit of being able to merge nodes, and our implementation of restricted decision diagrams does not \edit{currently} leverage equivalent \edit{node} merging.

\subsection{Relaxed Decision Diagrams}
\label{sec:background:relaxed}

There are many methods of constructing relaxed decision diagrams, and many heuristic decisions that must be made when doing so. In this paper, we focus on the method described by Cire \& van Hoeve (2013) for sequencing problems. As opposed to the top-down construction described in Section \ref{sec:background:restricted}, here $\underline{\mathcal{M}}$ will be constructed by separation. Constructing DDs by separation uses $\underline{\mathcal{M}}$ as a domain store over which constraints can be propagated. This method starts with a weak relaxation, and then strengthens it by splitting nodes until each layer is either exact, or has a width equal to $w_m$. The algorithm begins from a $1$-width MDD with an arc from the node on layer $\ell_i$ to the node on layer $\ell_{i+1}$ for each element that can be placed at position $\ell_i$ in the sequence. Thus, even though each layer only has one node, there can be several arcs between layers (see the relaxed diagram in Figure \ref{fig:basicDDs}). Then a node $u$ is selected and split to strengthen the relaxation. The process of splitting $u$ involves creating two new nodes $u^\prime_1$ and $u^\prime_2$, and then distributing the \emph{in} arcs of $u$ between $u^\prime_1$ and $u^\prime_2$. Then for each \emph{out} arc $a_{uv}$ from $u$, arcs $a_{u^\prime_1v}$ and $a_{u^\prime_2v}$ are added such that $a_{uv}$, $a_{u^\prime_1v}$ and $a_{u^\prime_2v}$ all have the same label. Finally $u^\prime_1$ and $u^\prime_2$ are filtered to remove infeasible and sub-optimal arcs. A collection of filtering rules are used to check each arc. As an example, given a feasible solution to $\mathcal{P}$ with objective value $z_{opt}$, an arc $a$ can be removed if all paths containing $a$ have an objective value greater than $z_{opt}$. The full process of identifying which arcs can be removed is detailed in \shortcite{MDDForSP}, and is not replicated here.

The following notation and definitions are critical to understanding these algorithms. Let $All_u^\downarrow$ be the set of arc labels that appear in every path from $r$ to $u$. Let $Some_u^\downarrow$ be the set of arc labels that appear in at least one path from $r$ to $u$. Let $All_u^\uparrow$ and $Some_u^\uparrow$ be defined as above, except that they refer to paths from $u$ to $t$. Let $\mathscr{J}$ be the set of all possible arc labels. For the SOP, we define an \emph{exact} node $u$ as a node where $Some_u^\downarrow = All_u^\downarrow$ and all arcs ending at $u$ originate from exact nodes. Intuitively, a node $u$ is exact if all paths to $u$ contain the same set of labels, and all parents of $u$ are exact. \edit{The $All$ and $Some$ notation was originally introduced by introduced by Andersen et al. (2007)}. Algorithm \ref{algo:RelaxedAlgo} formalizes the process of strengthening $\underline{\mathcal{M}}$.

\begin{algorithm}[h!t]
\SetAlgoLined
 Let $\underline{\mathcal{M}}$ be an MDD such that $Sol(\underline{\mathcal{M}}) \supseteq Sol(\mathcal{P})$\\
 \For{layer $L_j \in \underline{\mathcal{M}}$ from $j=1$ to $j=n$}{
 	 \While{$|L_j| < w_m$ and $\exists$ some node $ y \in L_j$ such that $y$ is not exact}{
 	 	$\mathscr{J} \leftarrow $ getAssignmentOrdering($\mathcal{P}$)\\
 	 	\small\emph{The getAssignmentOrdering() function returns a heuristically defined ordering of the values that can be assigned to decision variables}\normalsize\\
 	 	\For{$\phi \in \mathscr{J}$ \textbf{while} $|L_j| < w_m$}{
 			$S \leftarrow $ selectNodes($L_j$,$\phi$)\\
 			\small\emph{The selectNodes() function returns the set of nodes $u \in L_j$ such that $\phi \in Some_u^\downarrow \backslash All_u^\downarrow$}\normalsize\\
 			\For{$u \in S$ \textbf{while} $|L_j| < w_m$}{
 				Create two new nodes $u_1^\prime,u_2^\prime$\\
 				$L_j \leftarrow (L_j \cup \{u_1^\prime,u_2^\prime\})$\\
 				\ForEach{arc $a_{vu}$}{
					\eIf{$\phi \in (All_v^\downarrow ~ \cup$ the label of $a)$}
						{Redirect $a$ such that $a_{v u_1^\prime}$}
						{Redirect $a$ such that $a_{vu_2^\prime}$} 		
 				}
 				\ForEach{arc $a_{uv}$}{
 				    Create arcs $a_{u_1^\prime v}$ and $a_{u_2^\prime v}$ such that $label(a_{uv}) = label(a_{u_1^\prime v}) = label(a_{u_2^\prime v})$\\
 					filter($a_{u_1^\prime v}$), filter($a_{u_2^\prime v}$)\\
 					\small\emph{filter($a$) runs a list of quick checks to see if an arc can be removed}\normalsize\\
 				}
 				$L_j \leftarrow (L_j\backslash u)$\\
 			}
 		}
 	}
 }
 \Return{$\underline{\mathcal{M}}$}
 \caption{Refining Decision Diagrams for Sequencing~\shortcite{MDDForSP}}
  \label{algo:RelaxedAlgo}
\end{algorithm}

Deciding which nodes to split, and how to split them, are heuristic decisions with a significant impact on the bound that can be achieved without exceeding $w_m$~\shortcite{DDforO}. The algorithm discussed here selects nodes that can be split into equivalency classes, such that every path to the new node contains a certain label. Deciding which equivalency classes to produce first is another heuristic decision. The details of ordering the importance of the labels are specific to the problem being solved, and are not discussed here. However, it is important to note that the ordering for this implementation is static, and does not change between iterations. 

\subsection{Branch-and-Bound with Decision Diagrams}
\label{sec:background:bnb}

In a typical branch-and-bound algorithm, the branching takes place by splitting on the domain of the  variables. With decision diagrams, the branching takes place on the nodes themselves by selecting a set of exact nodes to represent the problem. The solver outlined by Bergman et al. (2016b) defines an exact node as a node $u$ for which every path from $r$ to $u$ ends in an equivalent state. As mentioned above, we can be more specific when applying this to sequencing problems, and define an exact node $u$ as a node where $Some_u^\downarrow = All_u^\downarrow$ and all arcs ending at $u$ originate from exact nodes. An \emph{exact cutset} is defined as a set of exact nodes that contain every path from $r$ to $t$. Let $\underline{\mathcal{M}}(u)$ be a relaxed decision diagram with root $u$, and let $\overline{\mathcal{M}}(u)$ be a restricted decision diagram with root $u$. The branch-and-bound algorithm for MDDs proceeds by selecting an exact cutset of $\underline{\mathcal{M}}$, and using each node $u$ in the cutset as the root for a new restricted decision diagram $\overline{\mathcal{M}}(u)$ and relaxed decision diagram $\underline{\mathcal{M}}(u)$. A node can be removed from the queue if the relaxation of that node is not better than the best known solution to $\mathcal{P}$, otherwise the exact cutset of the new node is added to the queue, and the process repeats until the queue is empty. This is detailed by Algorithm \ref{algo:BnB}.

Gillard et al. 2021) expanded on Algorithm \ref{algo:BnB} by incorporating a local search. A heuristic is used to quickly calculate a \emph{rough relaxed bound}\footnote{Gillard et al.~\shortcite{ddo} call the value \emph{rough upper bound}, but since we are testing a minimization problem in this paper, we use the term \emph{rough relaxed bound} instead.} at each node, and if the length of the shortest path to that node plus the rough relaxed bound is worse than the best known solution, the node can be removed. More formally, let $rrb(u)$ be a rough relaxed bound on $\mathcal{P}$ starting from node $u$, and let $z_{opt}$ be the value of best known solution so far. If $T^*(u)+rrb(u) > z_{opt}$, the node can be removed. They also provide evidence that if $rrb(u)$ is inexpensive to compute, it can be used to filter nodes in $\overline{\mathcal{M}}$ and $\underline{\mathcal{M}}$. The method of using a rough relaxed bound to trim nodes is used in this paper, but the details are problem specific and are discussed in a later section. 

\begin{algorithm}[!ht]
\SetAlgoLined
Let $\underline{\mathcal{M}_{uu^\prime}}$ be a partial diagram with root $u$ and terminal $u^\prime$\\
Let $v^*(u)$ be the lower bound of $\mathcal{P}$ resulting from starting at node $u$\\
Let $z_{opt}$ be the value of the best known solution\\
 $Q = \{r\}$ \\
 $v^*(r) \leftarrow 0$\\
 $z_{opt} \leftarrow \infty$\\
 \While{$Q \neq \emptyset$}{
 	$u \leftarrow$selectNode($Q$),  $Q \leftarrow Q\backslash\{u\}$\\
 	$\overline{\mathcal{M}} \leftarrow \overline{\mathcal{M}}(u)$\\
 	 \If{$v^*(\overline{\mathcal{M}}) < z_{opt}$} {$z_{opt} \leftarrow v^*(\overline{\mathcal{M}})$\\}
 	 \If{$\overline{\mathcal{M}}$ is not exact}{
 	 	$\underline{\mathcal{M}} \leftarrow \underline{\mathcal{M}}(u)$\\
 	 	\If{$v^*(\underline{\mathcal{M}}) < z_{opt}$} {
 	 		$S \leftarrow $ exactCutset($\underline{\mathcal{M}}$)\\
 	 		\ForEach{$u^\prime \in S$}{
				let $v^*(u^\prime) = v^*(u) + v^*(\underline{\mathcal{M}_{uu^\prime}})$\\
				$Q \leftarrow Q \cup u^\prime$\\
 	 		}
 	 	}
 	 }
 }
 \Return{$z_{opt}$}
 \caption{Decision Diagram based Branch-and-Bound  (BnB)~\shortcite{BnB2}}
 \label{algo:BnB}
\end{algorithm}

\section{Peel-and-Bound Algorithm}
\label{sec:pnb}

The motivation for peel-and-bound stems from an observation about Algorithm~\ref{algo:RelaxedAlgo}. When implemented in a branch-and-bound structure, a large portion of the work done while generating each $\underline{\mathcal{M}}$ is repeated at every iteration. Creating the relaxation for some exact node $u$ in the queue requires creating a $1$-width decision diagram, iterating over each layer from the top down, and splitting nodes in a predetermined order. The static order of node splits means that for each node $y$ such that $\ell_y > \ell_u$, the first equivalency class created when splitting $y$ is the same in $\underline{\mathcal{M}}(r)$ and $\underline{\mathcal{M}}(u)$. The existing arcs for both diagrams will be sorted in the same way, and the only difference is the possibility of filtering arcs in $\underline{\mathcal{M}}(u)$ that could not be filtered in $\underline{\mathcal{M}}(r)$ due to the added constraint that all paths must pass through $u$. The extra filtered arcs are the reason that $\underline{\mathcal{M}}(u)$ may produce a stronger bound than $\underline{\mathcal{M}}(r)$. However, because equivalency classes are chosen in the same order each time, many arcs that were filtered while constructing $\underline{\mathcal{M}}(r)$ will also be filtered again while constructing $\underline{\mathcal{M}}(u)$. There is a sub-graph of $\underline{\mathcal{M}}(r)$, induced by node $u$, that contains all of the paths that will be encoded in $\underline{\mathcal{M}}(u)$, but does not contain the arcs that are filtered from both diagrams during construction. Thus, less work needs to be performed at each iteration of branch-and-bound by starting from that sub-graph instead of a $1$-width diagram. If the split order is static, the same diagram is generated starting from either the $1$-width diagram, or the sub-graph induced by $u$. If the split order changes between branch-and-bound iterations, the sub-graph induced by $u$ is still a valid relaxation, but the generated diagram will differ from one that began at width $1$. 

\begin{figure}[!ht]
    \centering
    \begingroup
        \setlength{\tabcolsep}{10pt}
        \begin{tabular}{c c c c}
            \includegraphics[scale=.67]{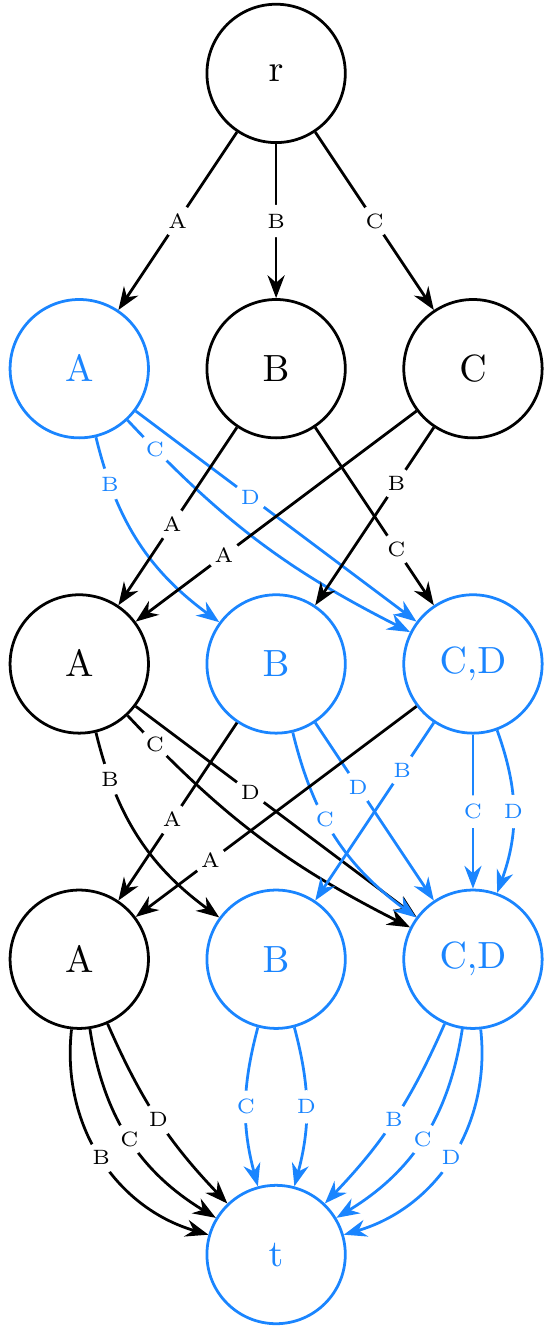} & \includegraphics[scale=.67]{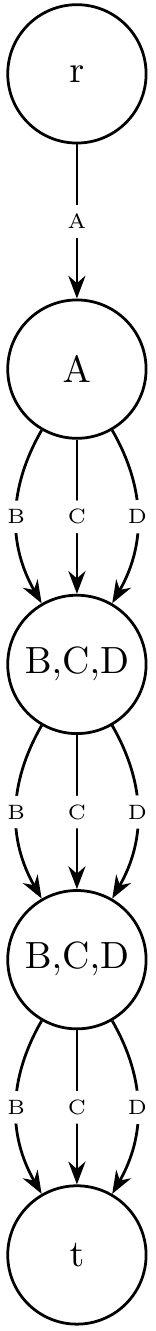} & \includegraphics[scale=.67]{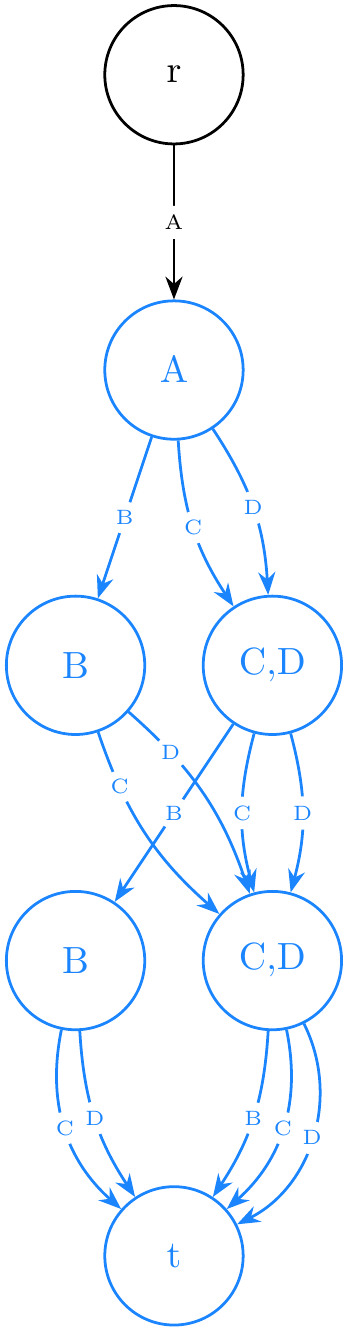} & \includegraphics[scale=.67]{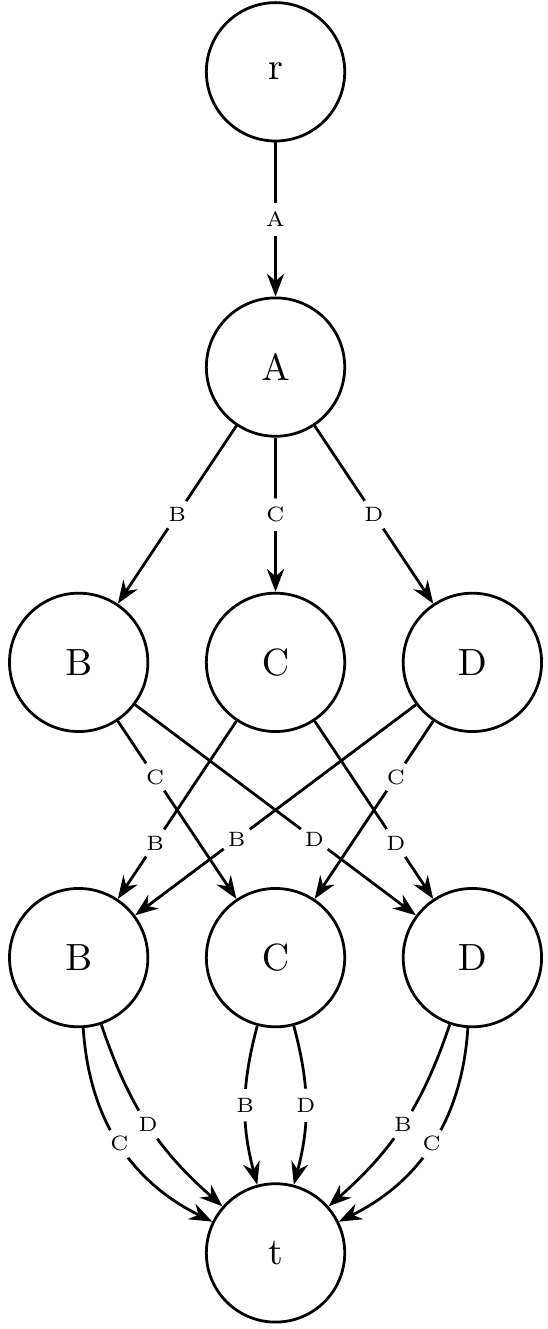}\\
            $\underline{\mathcal{M}}(r)$ & $\underline{\mathcal{M}}(A)$ Stage 1 & $\underline{\mathcal{M}}(A)$ Stage 2 & $\underline{\mathcal{M}}(A)$ Stage 3
        \end{tabular}
    \endgroup
    \caption{Example of \edit{a sub-graph (shown in blue), and the associated relaxed decision diagram with the same root for a SOP instance. The objective is to order the elements $[A,B,C,D]$, subject to the precedence constraint that $A$ must precede $D$, an alphabetical ordering heuristic, and $w_m = 3$.} }
    \label{fig:SubDiagram}
\end{figure}

Consider a SOP instance where the goal is to order the elements $[A,B,C,D]$, subject to the precedence constraint that $A$ must precede $D$, an alphabetical ordering heuristic, and $w_m = 3$. Figure~\ref{fig:SubDiagram} shows $\underline{\mathcal{M}}(r)$, and $\underline{\mathcal{M}}(A)$ in three stages. The first stage is the initial $1$-width diagram. The second stage is after one split on each layer, and the third stage is the complete diagram. The sub-graph shared by $\underline{\mathcal{M}}(r)$ and $\underline{\mathcal{M}}(A)$ is highlighted in blue, indicating that in this case the first two splits could have been read from $\underline{\mathcal{M}}(r)$ instead of being re-created from scratch.

\begin{figure}[!ht]
    \centering
    \begingroup
        \setlength{\tabcolsep}{20pt}
        \begin{tabular}{c c}
            \includegraphics[scale=.59]{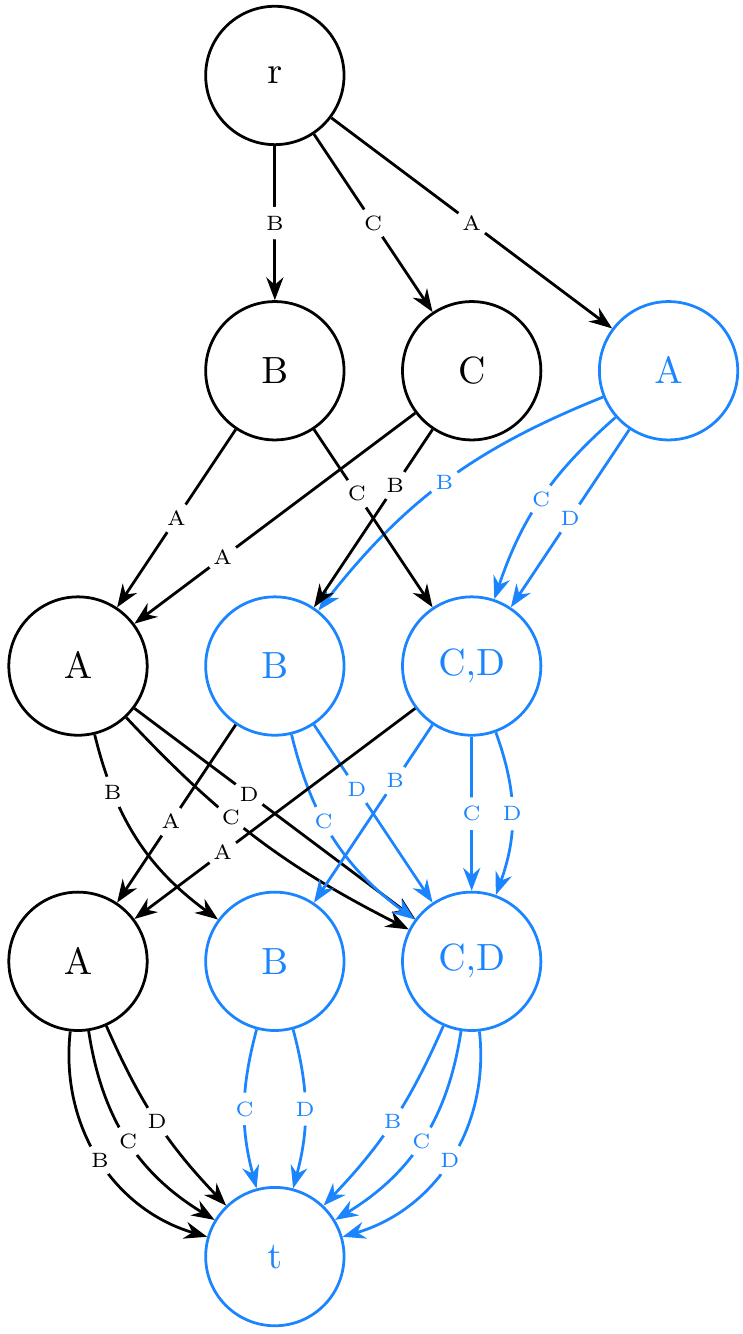} & \includegraphics[scale=.59]{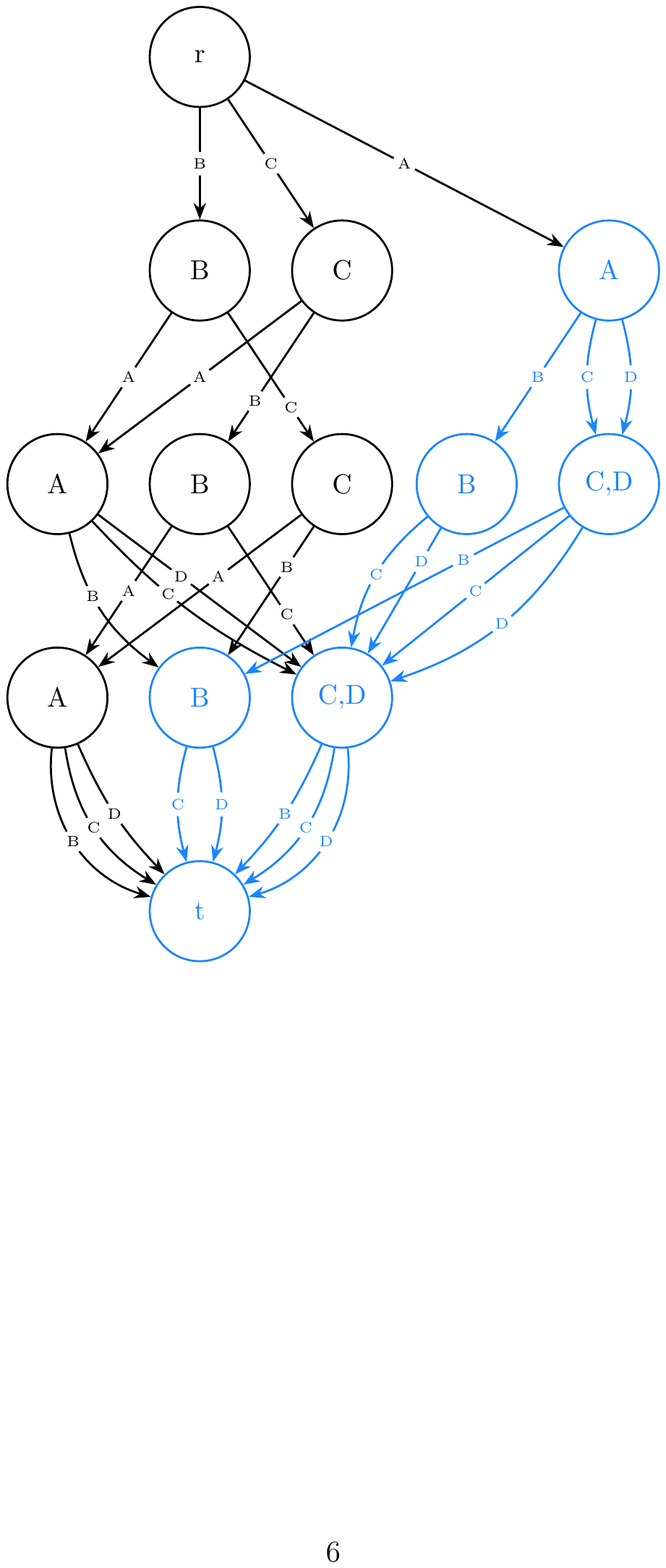}\\
            (1) & (2) \\
             & \\
            \includegraphics[scale=.59]{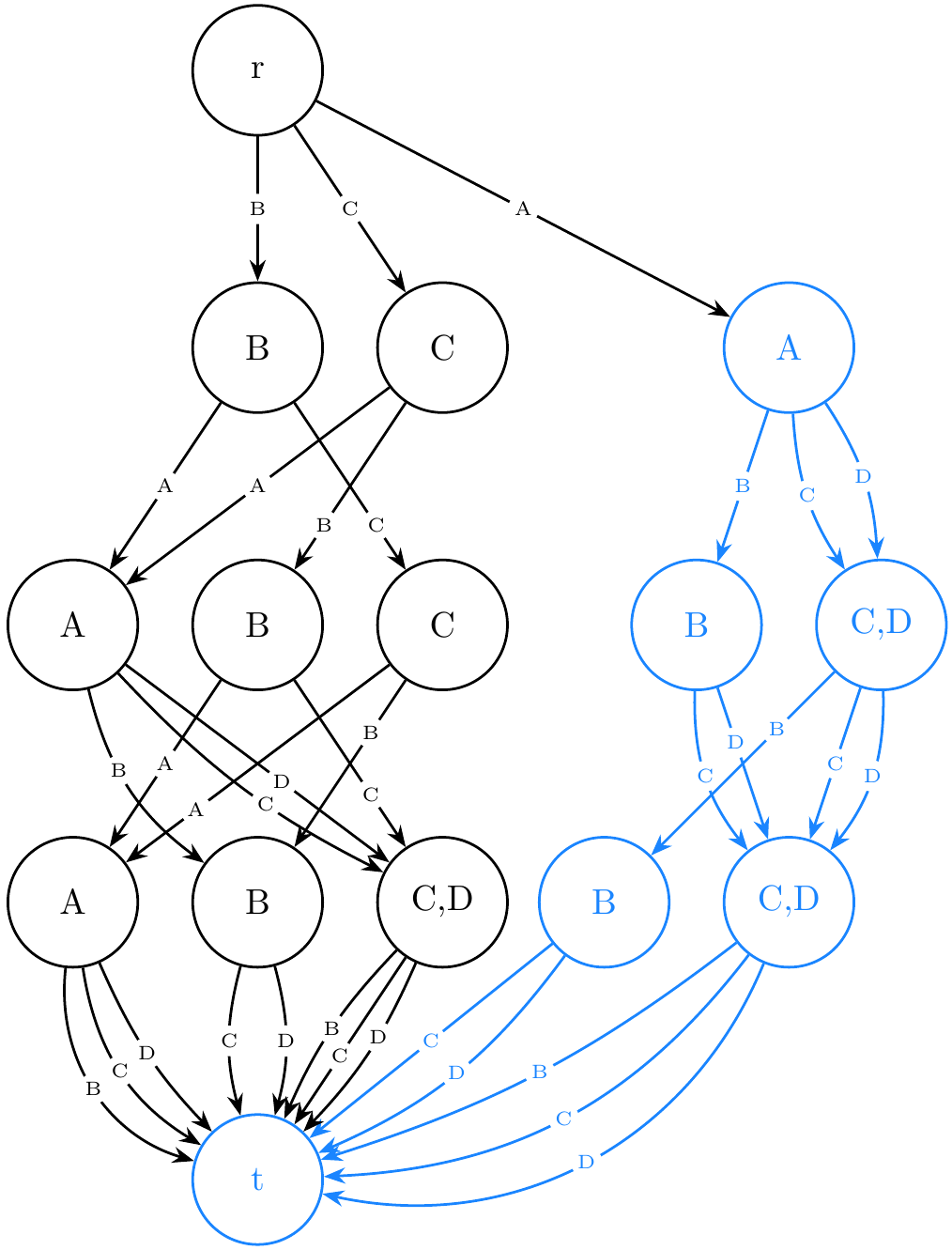} & \includegraphics[scale=.59]{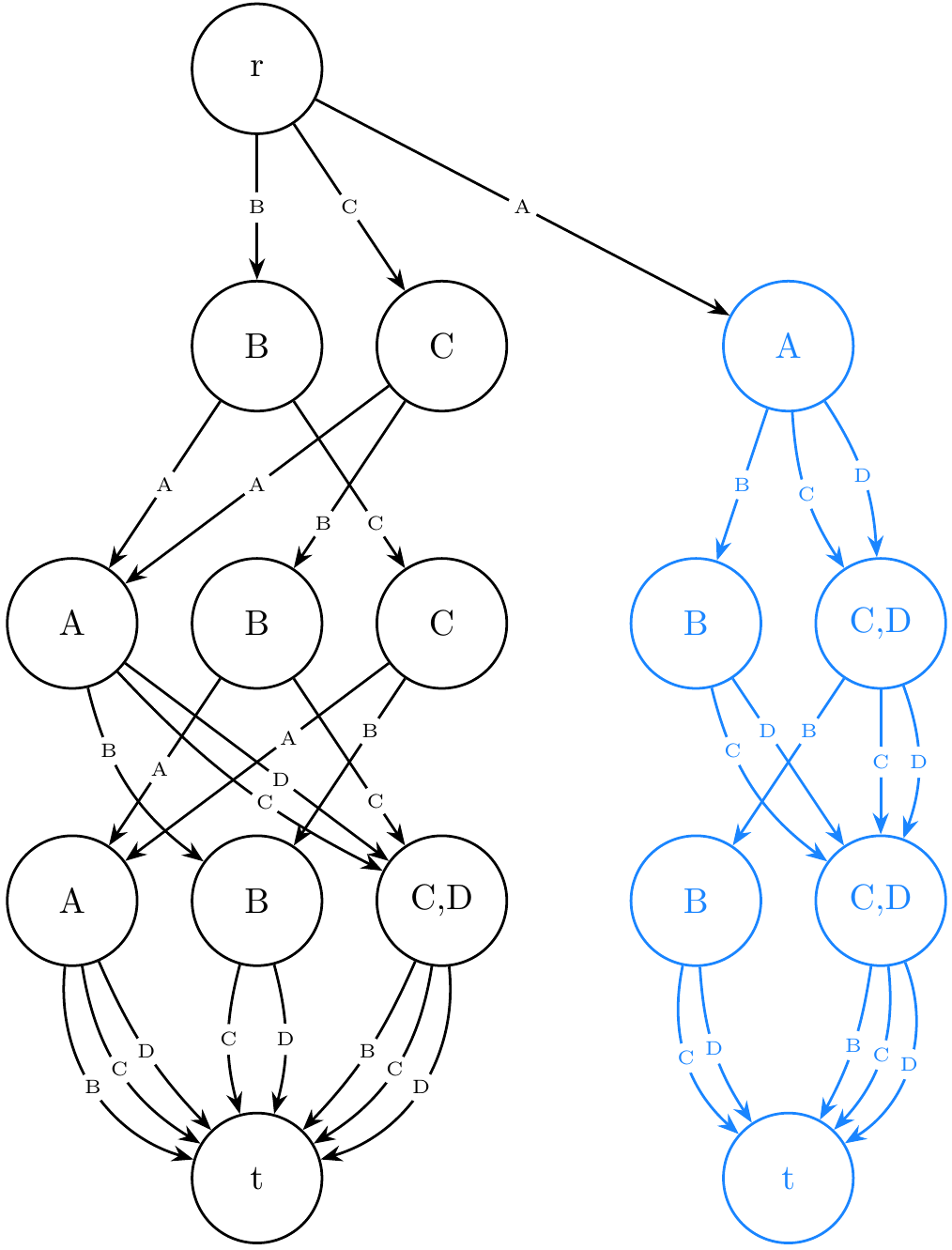} \\
            (3) & (4)\\
        \end{tabular}
    \endgroup
    \caption{An example of a peel operation. In (1), $A$ is selected to induce the peel process and removed from the the original diagram ($\underline{\mathcal{M}}(r)$ from Figure \ref{fig:SubDiagram}). In (2) the arcs that connect $A$ to the original diagram are moved to copies of the nodes they originally ended at, and infeasible arcs are filtered. In (3) and (4) the process is repeated until the diagrams are disconnected.}
    \label{fig:PnB}
\end{figure}

This mechanism can be embedded into a slightly modified version of the standard branch-and-bound algorithm based on decision diagrams (Algorithm~\ref{algo:BnB}). In peel-and-bound, the queue stores diagrams instead of nodes. After the initial relaxation $\underline{\mathcal{M}}(r)$ is generated, the entire diagram is placed into the queue $Q$ such that $Q = \{\underline{\mathcal{M}}(r)\}$. Then a diagram $\underline{\mathcal{M}}(u)$ is selected from $Q$ (for the first iteration $\underline{\mathcal{M}}(u) = \underline{\mathcal{M}}(r)$). However, instead of selecting an exact cutset of $\underline{\mathcal{M}}(u)$, a single exact node $e$ from $\underline{\mathcal{M}}(u)$ is selected. The process of selecting a diagram and exact node are heuristic decisions that are discussed in Section~\ref{sec:pnb:implementation}. The process of \emph{peeling} $e$ is as follows. Create an empty diagram $\underline{u}$, remove $e$ from $\underline{\mathcal{M}}(u)$, and then put $e$ into $\underline{u}$ such that $e$ is the root of $\underline{u}$, and the arcs leaving $e$ still end in $\underline{\mathcal{M}}(u)$. Then for each node $y$ in $\underline{\mathcal{M}}(u)$ with an \emph{in} arc that originates in $\underline{u}$, a new node $y^\prime$ is made and added to $\underline{u}$. Each \emph{in} arc $a_{oy}$ of $y$ that originates in $\underline{u}$ is removed and then arc $a_{oy^\prime}$ is added to $\underline{u}$. Then the \emph{out} arcs of $y$ and $y^\prime$ are filtered using the same \emph{filter} function as Algorithm~\ref{algo:RelaxedAlgo}. The process of removing and adding arcs is repeated until there are no arcs ending in $\underline{\mathcal{M}}(u)$ that originate in $\underline{u}$. This procedure accomplishes a top-down reading of the sub-graph induced by $e$, and potentially strengthens $\underline{\mathcal{M}}(u)$ by removing nodes and arcs in the process. If the shortest path through the modified $\underline{\mathcal{M}}(u)$ is less than the best known solution, $\underline{\mathcal{M}}(u)$ is put back into $Q$. The diagram $\underline{u}$ is relaxed using Algorithm \ref{algo:RelaxedAlgo}. Let $\underline{\mathcal{M}}(\underline{u})$ be the result; then if the shortest path through the refined diagram $\underline{\mathcal{M}}(\underline{u})$ is less than the best known solution, $\underline{\mathcal{M}}(\underline{u})$ is added to $Q$. The whole procedure is repeated until there are no nodes left in the queue ($Q = \emptyset$). A peel operation is illustrated and explained in Figure \ref{fig:PnB}. Peel-and-bound is formalized in 
Algorithm~\ref{algo:PnB}, and the peel operation is formalized in Algorithm~\ref{algo:Peel}.

\begin{algorithm}[!ht]
\SetAlgoLined
Let $v^*(u)$ be the lower bound of $\mathcal{P}$ resulting from starting at node $u$\\
Let $z_{opt}$ be the value of the best known solution\\
 $Q = \{\underline{\mathcal{M}}(r)\}$ \\
 $z_{opt} \leftarrow \infty$\\
 \While{$Q \neq \emptyset$}{
    $\mathscr{D} \leftarrow$selectDiagram($Q$), $Q \leftarrow Q\backslash\{\mathscr{D}\}$\\
    $u \leftarrow$selectExactNode($\mathscr{D}$)\\
    $\underline{u},\mathscr{D}^*\leftarrow$ peel($\mathscr{D}$, $u$) \emph{(See Algorithm \ref{algo:Peel})}\\
    \If{$v^*(\mathscr{D}^*) < z_{opt}$}{
        $Q \leftarrow Q \cup \{\mathscr{D}^*\}$
    }
    $\overline{\mathcal{M}} \leftarrow \overline{\mathcal{M}}(u)$\\
    \If{$v^*(\overline{\mathcal{M}}) < z_{opt}$} {$z_{opt} \leftarrow v^*(\overline{\mathcal{M}})$\\}
 	 \If{$\overline{\mathcal{M}}$ is not exact}{
 	 	$\underline{\mathcal{M}} \leftarrow \underline{\mathcal{M}}(\underline{u})$\\
 	 	\If{$v^*(\underline{\mathcal{M}}) < z_{opt}$} {
 	 		$Q \leftarrow Q \cup \{\underline{\mathcal{M}}\}$
 	 	}
 	 }
 }
 \Return{$z_{opt}$}
 \caption{Peel-and-Bound (PnB) Algorithm}
 \label{algo:PnB}
\end{algorithm}

\begin{algorithm}[!ht]
\SetAlgoLined
    Let $in(u)$ for some node $u$ be the set of arcs that end at node $u$\\
    Let $out(u)$ for some node $u$ be the set of arcs that originate from node $u$\\
    Let $in(\mathcal{M})$ for some MDD $\mathcal{M}$ be the set of arcs that end in $\mathcal{M}$\\
    Let $out(\mathcal{M})$ for some MDD $\mathcal{M}$ be the set of arcs that originate in $\mathcal{M}$\\
    \textbf{input:} a relaxed MDD $\mathscr{D}$, and an exact node $u$ in $\mathscr{D}$\\
    \textbf{output:} a relaxed MDD $\underline{u}$ peeled from $\mathscr{D}$, and what remains of relaxed MDD $\mathscr{D}$\\
    Let $\underline{u}$ be an empty decision diagram\\
    $in(u) \leftarrow \emptyset$\\
    $\mathscr{D} \leftarrow \mathscr{D}\backslash u$\\
    $\underline{u} \leftarrow u$\\
    \While{$in(\mathscr{D})\cap out(\underline{u}) \neq \emptyset$}{
        \ForEach{node $m \in \mathscr{D}$ with an in arc that originates in $\underline{u}$}{
            create a new node $m^\prime$, and add it to $\underline{u}$\\
            \ForEach{arc $a \in in(m)$ that originates in $\underline{u}$}{
                change the destination of $a$ to $m^\prime$\\
                filter($a$)
            } 
            \ForEach{arc $a \in out(m)$}{
                filter($a$)
            }
        }
    }
    \While{$\exists$ some node $m \in D$ with $in(m)=\emptyset$ or $out(m) = \emptyset$ (excluding $r$ and $t$)}{
            $in(m) \leftarrow \emptyset$\\
            $out(m) \leftarrow \emptyset$\\
            $\mathscr{D} \leftarrow \mathscr{D} \backslash \{m\}$
    }
    \Return{($\underline{u}$, $\mathscr{D}$)}
    \caption{Peeling process used in Algorithm \ref{algo:PnB}}
    \label{algo:Peel}
\end{algorithm}

\subsection{Complexity Analysis}
\label{sec:pnb:complexity}

Separating each node $u$ during a peel requires creating a new node $u^\prime$, moving the \emph{in} arcs of $u$ that originate in the peeled diagram $\underline{u}$ to \emph{$u^\prime$}, copying the \emph{out} arcs of $u$ to $u^\prime$, and then filtering the \emph{out} arcs of $u$ and $u^\prime$. Creating a new node in our implementation has a time in $\mathcal{O}(n)$, where $n$ is the number of elements in the sequence (or the number of possible decisions in a more generic problem), due to storing state information that has a size in $\mathcal{O}(n)$ (such as $All_{u^\prime}^\downarrow$). However, it is possible that in other applications the size of a node is in $\mathcal{O}(1)$. The number of \emph{in} arcs of $u$ is at most $w$, although this worst case is unlikely in practice because it requires $\underline{u}$ to have width $w$ and for each node in $\underline{u}$ on layer $\ell_{u}-1$ to have an arc ending at $u$. Thus, moving the \emph{in} arcs of $u$ has a time in $\mathcal{O}(w)$. The number of \emph{out} arcs of $u$ is at most $n$, and each arc has a size in $\mathcal{O}(1)$, so copying the \emph{out} arcs has a time in $\mathcal{O}(n)$. Each individual filtering process has a time in $\mathcal{O}(1)$ as it uses only existing state information from $u$ and $u^\prime$, and it is performed on the at most $2n$ \emph{out} arcs of $u$ and $u^\prime$. Thus, filtering the \emph{out} arcs has a time in $\mathcal{O}(n)$. Therefore, separating one node during the peel process has a time in $\mathcal{O}(n+w)$. Separations during a standard relaxation procedure require selecting a node ($\mathcal{O}(w)$), making a new node ($\mathcal{O}(n)$), partitioning the \emph{in} arcs ($\mathcal{O}(nw)$), copying the \emph{out} arcs ($\mathcal{O}(n)$), and filtering the \emph{out} arcs ($\mathcal{O}(n)$). The reason that there can be more \emph{in} arcs during a standard relaxation procedure is because the nodes in a $1$-width diagram can have \emph{in} arcs with different labels coming from the same node, whereas the structure of the diagram during a peel guarantees that each node $u$ can have only one \emph{in} arc from each node on the layer $\ell_{u}-1$. Thus, the total time for a separation in a standard relaxation is in $\mathcal{O}(nw)$.

The maximum number of separations during a peel is the maximum number of nodes in the peeled diagram. A peeled diagram can have at most $(n-3) \times w + 2$ nodes, and thus the number of nodes is in $\mathcal{O}(nw)$. Therefore, the entire peel process has a time in $\mathcal{O}(n^2w + nw^2)$. The maximum number of separations during a standard relaxation is the exact same as during a peel, since the resulting diagram will be the same size. Thus, the standard relaxation has a total time in $\mathcal{O}(n^2w^2)$. However, peel-and-bound uses a peel to generate some fraction of the nodes, then a standard relaxation to generate the rest. Let $\alpha$ be the percent of nodes that are peeled during the peel. It follows that the total time for an iteration of peel-and-bound is in $\mathcal{O}(\alpha(n^2w + nw^2)+(1-\alpha)(n^2w^2))$. Therefore, the larger that $\alpha$ grows, the more time peel-and-bound saves over branch-and-bound.

\subsection{Advantages and Implementation Decisions}
\label{sec:pnb:implementation}
\subsubsection{Node Selection}
\label{sec:pnb:implementation:nodeselection}
The branch-and-bound algorithm proposed by Bergman et al. (2016b) requires selecting an exact cutset of $\underline{\mathcal{M}}$. Peel-and-bound requires selecting a diagram from the queue, and an exact node to start the peel process. The choice of node has a substantial impact on how quickly the process converges to an optimal solution, because it serves two purposes simultaneously. As discussed earlier, the first purpose of peeling is to avoid recreating a portion of the diagram at each iteration. The second purpose is to strengthen the overall relaxation. Let $\underline{u}$ be a diagram peeled from $\underline{\mathcal{M}}$, and let $\underline{\mathcal{M}}^*$ be $\underline{\mathcal{M}}$ after the peel operation. If $Sol(\mathcal{P}) \subseteq Sol(\underline{\mathcal{M}})$ then $Sol(\mathcal{P}) \subseteq Sol(\underline{\mathcal{M}}^*)\cup Sol(\underline{u})$. The only step of peel-and-bound that removes paths is the \emph{filter} step, which only removes an arc if no feasible solutions can pass through that arc. If the node the peel is induced from contains the shortest path through $\underline{\mathcal{M}}$, then there will be a new shortest path through $\underline{\mathcal{M}}^*$ with $T^*(\underline{\mathcal{M}}^*) \geq T^*(\underline{\mathcal{M}})$. Similarly after peeling, the peeled diagram is going to be strengthened and $T^*(\mathcal{M}(\underline{u})) \geq T^*(\underline{u})$. Therefore, when implementing the \emph{selectDiagram} and \emph{selectExactNode} functions from Algorithm \ref{algo:PnB}, we propose selecting the diagram $\mathscr{D}$ containing the most relaxed bound, and an exact node from $\mathscr{D}$ that contains $z^*(\mathscr{D})$ at each iteration. Using these parameters, the peel step of peel-and-bound strengthens the relaxed bound of $\mathcal{P}$, in addition to providing a stronger initial diagram to use when generating $\mathcal{M}(\underline{u})$.

We originally proposed two heuristics for selecting a node from $\mathscr{D}$ that contains $z^*(\mathscr{D})$. The first heuristic picks the first node in the shortest path through the diagram with at least one child that is not exact, we call this the \emph{last exact node}. The second heuristic picks the \emph{frontier} node, the highest-index exact node that contains $z^*(\mathscr{D})$. Taking the last exact node is more of a breadth-first search that peels a large set of nodes where all of the paths to the nodes in the set share a beginning with $z^*(\mathscr{D})$. In contrast, taking the frontier node is more of a depth-first search, taking fewer nodes and exploring those nodes at greater depth. We have since experimented with a third heuristic where we pick what we call the \emph{maximal} node. The maximal node is simply the node on the second layer that contains $z^*(\mathscr{D})$. This peels as many nodes as possible while still picking a node that contains $z^*(\mathscr{D})$.

\subsubsection{Handling Non-Separable Objective Function}
\label{sec:pnb:implementation:unlabeled}

\edit{When a problem is represented with an Integer Program, the value or cost of making a decision is separate for each variable; in other words the objective function is separable. However, when using MDDs to represent a problem, the objective function can be non-separable. In the example SOP used in Section \ref{sec:background}, the cost of an arc leaving a node $u$ is dependent on the labels of the arcs that end at $u$. If the arcs have conflicting labels, then arcs leaving $u$ may not have an exact cost.} Cire \& van Hoeve (2013) proposed that each iteration of Algorithm \ref{algo:RelaxedAlgo} starts from a $1$-width MDD. However, for peel-and-bound with a non-separable objective function, starting from a $1$-width MDD poses a problem. The arcs in such a diagram do not have exact values, because they are dependent on the state of the node they originate from. As nodes are peeled, the values of those arcs must be updated, and the operation becomes computationally expensive at scale. 

\edit{When solving sequencing problems where the non-separability of transition costs arises from their dependency on the previous element in the sequence, }this problem can be avoided by creating the initial diagram using a structure where all of the arcs ending at a given node have the same label. The resulting initial diagram has a width of $n$, and each node on the layer is assigned to one state $s \in \{1,...,n\}$. Then every possible feasible arc between consecutive layers is added. Thus, the nodes of $\underline{\mathcal{M}}$ do not have relaxed states, and each arc can only take one possible value. Starting from such a diagram not only removes the need to update arc values, it ensures that every arc generated during peel-and-bound is an exact copy of an arc that exists in the initial diagram, since arcs are only copied or removed, never updated or added. Using this structure makes implementation easier by removing the need to store any information on the arcs at all. As each node will only have one state $s$, the label and weight of an arc is implied by the state of the node it originates from, and the node it points to. This allows all information in a diagram to be stored on the nodes, and thus only the information on the nodes ever needs to be read or updated. 

\edit{When transition costs are dependent on factors other than the previous element in a sequence, such as the Time-Dependent TSP (TSP-TD), modifications to the described method become necessary. In general an exact cost may not be calculable, but a relaxed bound on cost can be made calculable by storing the relevant values as a range at each node. Consider an MDD for the TSP-TD where nodes $u$ and $v$ are connected by an arc. Let $a_u$ to $b_u$ be the range of possible times one can arrive at $u$. While an exact cost for the arc from $u$ to $v$ is not available, a relaxed bound can be found by varying the time from $a_u$ to $b_u$, and then using the best value found as the cost. As the MDD is peeled, the gap between $a_u$ and $b_u$ may grow smaller, but it will never grow larger. When the gap decreases, the previously calculated cost is still a valid relaxed bound. Thus, the cost can be updated to retrieve an improved bound, but it is not necessary to do so with every peel operation. Many problems will require modifying this concept to fit their particular constraints, but this method is easy to adapt in most cases.} An alternative method of handling non-separable objective functions is explored by \shortcite{canonicalarcs,jobsequencing,jobsequencing2}.

\subsubsection{Parallelization}
\label{sec:pnb:implementation:parallelization}

The decision diagram based branch-and-bound shown in Algorithm \ref{algo:BnB} is particularly amenable to parallelization as shown by Bergman et al. (2014\edit{b}) and again by Gillard (2022). Algorithms seeking to parallelize must overcome the data-race problem. In other words, if multiple processors are working on a problem simultaneously, then there must be a process in place to stop them from trying to write to the same place in memory at the same time. For many algorithms, this poses a substantial challenge or creates substantial overhead. However, for both decision diagram based branch-and-bound, and peel-and-bound, the solution is simple. As a problem is being solved, nodes (or diagrams in the case of peel-and-bound) are placed into a processing queue. Each node/diagram represents a discrete problem that needs to be solved, and can be processed separately. Given access to a sufficiently large number of processors, each node/diagram added to the queue could be immediately dispatched to an available processor for processing. The only communication required between the processors is the current value of the best known solution. In theory, this could result in a linear improvement in time spent solving a problem when increasing the number of processors available, because $k$ processors may be able to process $k$ nodes, in the time it takes $1$ processor to process $1$ node. \edit{This process could also result in a superlinear speedup due to the non-deterministic processing order of elements in the queue. In the parallel implementation, one processor may find an improved bound which can then be utilized by the other processors. Therefore, elements from the queue could be processed using bounds that are better than the bounds used by the single-thread deterministic implementation. The ability to identify better bounds earlier on during parallel processing when processing an element could yield less elements that need to be processed overall.
} In practice, however, there are heuristic decisions that must be made that can have a large impact on solve time. 

The dilemma one encounters in implementation is the existence of a critical path in the solution finding process. To demonstrate this with an extreme example, consider a problem that does the following when solved using peel-and-bound on a single processor. Each time a node is peeled, one of the two resulting diagrams is solved and closed without requiring any additional peel operations. Then the single remaining diagram is processed again, a node is peeled, and the whole process repeats itself $m$ times. In this example, there are roughly $2m$ diagrams that need to be processed, but only $2$ are ever available at the same time. So in this case, $k$ processors would take exactly as long as $2$ processors to solve the problem. We propose two methods of handling this dilemma. The first is simple; reduce the maximum width of the diagrams. When solving a problem that is encountering this critical path problem, the work can be divided more equitably among the available processors by reducing the amount of work done at each iteration. When using a single processor, a higher maximum width is almost always more desirable as long as the decision diagrams can still be generated quickly, because the extra space makes it more likely that the diagram will be solved instead of producing more diagrams to add to the queue. Reducing the maximum width will increase the number of diagrams that need to be processed to solve the problem, but also produces those diagrams more quickly. A width too low can generate an enormous number of diagrams without closing any of them. The optimal width to use is one that generates enough diagrams for all available processors to have consistent work, but does not generate a large backlog of work. 

The second method of handling the critical path issue is to use the peel process to redistribute work as needed. Each time a processor is available and not being used, the peel procedure can peel off a sub-diagram for that processor to work on. The downside of this is that often a lot of the work that only needed to be performed once, will occur on both diagrams. The first method of simply lowering the maximim width accomplishes the same goal, but each diagram reaches its assigned maximum width before it is further processed, so arcs that can be processed out, only need to be processed out of one diagram. Any implementation of this second method would likely require more heuristic decisions to ensure the task scheduler distributes the work in a useful way.

\subsubsection{Embedded Restricted Decision Diagrams}
\label{sec:pnb:implementation:betterdds}
In both decision diagram based branch-and-bound (Algorithm \ref{algo:BnB}), and peel-and-bound (Algorithm \ref{algo:PnB}), at each branch a restricted decision diagram is created before the relaxed decision diagram. This is useful not just for solution finding, but also because the process of creating a restricted decision diagram is often several orders of magnitude faster than the process of refining a relaxed decision diagram. When the restricted decision diagram is exact, the relaxed decision diagram can be closed without any additional processing. Peel-and-bound provides an opportunity to leverage relaxed decision diagrams to improve the associated restricted decision diagrams. 

The following idea extends the work done in Copp\'e et al. (2022). \edit{Each path from the root to a node} in a restricted decision diagram represents a partial solution to the problem being solved. Any partial solution to the problem represents a partial path through a matching relaxed decision diagram. \edit{So each path to a node} $\overline{u}$ in a restricted decision diagram  $\overline{\mathcal{M}}$ can be mapped to exactly one node $\underline{u}$ in a relaxed decision diagram  $\underline{\mathcal{M}}$. When generating $\overline{\mathcal{M}}$ from the top-down, each node in $\overline{\mathcal{M}}$ creates a child node on the next layer for every element in its domain (before the layer is trimmed down to the maximum width). We incorporate multiple methods from the literature for trimming this domain, such as the rough relaxed bound we will discuss in Section \ref{sec:pnb:rrb}, but add our own here. Let $d(\overline{u})$ be the domain of $\overline{u}$; we set $d(\overline{u}) = d(\overline{u}) \cap d(\underline{u})$ before generating the child nodes of $\overline{u}$. There are two clear benefits to this strategy, that pair to the two reasons something can be in $d(\overline{u})$ but not $d(\underline{u})$. If an arc has been proven to be infeasible or sub-optimal, it will not be in $\underline{\mathcal{M}}$. However, the methods of that proof may not be available to nodes in $d(\overline{u})$. So our method prevents $\overline{M}$ from exploring nodes that have already been proven useless by $\underline{M}$. The other reason something can be in $d(\overline{u})$ but not $d(\underline{u})$ is that it was removed during a peel procedure performed on $d(\underline{M})$. Without using this intersection operation as a way to trim the domain, $\overline{\mathcal{M}}$ will search the entire solution space that starts from the same root as $\underline{\mathcal{M}}$, even if the peeled diagrams have already been closed as sub-optimal. Using this method, each restricted decision diagram will only search the \edit{solutions} encoded within the matching relaxed decision diagram. This means that each restricted decision diagram has a significantly improved chance of closing the relaxed decision diagram it maps to, because the solution space it must explore becomes smaller with each peel. 

\subsubsection{Search Diversification}
\label{sec:pnb:implementation:variety}
The structure of peel-and-bound yields another method of searching for solutions that forces increased diversification. The embedded restricted decision diagram described in Section \ref{sec:pnb:implementation:betterdds} takes advantage of the reduced search space embedded in the relaxed decision diagram, but it makes no effort to explore substantially different regions, and thus it is at risk of getting stuck in a local optimum. Here we propose a simple method of diversifying the solutions explored. Starting from the root, and moving down layer by layer, map each node $u$ in the relaxed decision diagram to the best feasible path that ends at $u$, and is a continuation of a path a parent of $u$ maps to. The obvious drawback is that many paths will becomes infeasible or sub-optimal, and many of the relaxed nodes might not map to a feasible path using this method, simply because the paths that were being explored in their parents were bad paths. To fix this, the number of paths stored can be expanded. Let $k$ be any positive integer; if each node maps to the $k$ best paths to that node, then as the value of $k$ increases there is a much higher likelihood of new, and better, solutions being found. However, as the number of nodes in the diagram can be quite large, even small values of $k$ can be computationally expensive. So this method can be a powerful tool for diversification, but it comes with a significant drawback in terms of compute time. It has the potential to be valuable if used just once at the beginning of the peel-and-bound process to search for initial solutions, but is unlikely to be useful if repeated often. This idea may also benefit from being combined with the large neighborhood search using restricted decision diagrams proposed by Gillard \& Schaus, (2022).

\subsection{Limitations and Handling Memory}
\label{sec:pnb:limitations}
The focus of this paper is sequencing problems, but peel-and-bound can be easily applied to other optimization problems. However, some existing MDD based methods conflict with peel-and-bound. For example, some MDD algorithms use a dynamic variable order~\shortcite{karahalios}, such that the variables the layers on $\overline{M}$ are mapped to in one iteration of branch-and-bound, are different in the next. Peel-and-bound as it is presented in this paper cannot be paired with a dynamic variable order. Furthermore, the method in this paper is specific to decision diagrams generated using separation. We believe the method can be extended to decision diagrams that use a merge operator, but it has not been shown here.

Memory limitations present a problem for peel-and-bound in theory, but not in practice. Each open diagram remains in the queue, and thus must be stored in memory. However, this problem can be handled in many ways; two are given here. A dynamic method of handling the problem is to start targeting large diagrams with bounds close to $z_{opt}$ as memory limitations start to become a problem. Such diagrams can usually be closed quickly, and subsequently removed from memory, freeing up space for the algorithm to continue. Alternatively, the diagrams with bounds closest to $z_{opt}$ can be deleted in favor of storing just the root node, then when they need to be processed, initial diagrams are generated for those once again. This method essentially falls back to branch-and-bound until memory limitations cease to be a problem. Additional approaches for working with memory limitations, and evidence that the problem can be handled efficiently, are presented by Perez \&
R\'egin, (2018).

\subsection{Integrating Rough Relaxed Bounds}
\label{sec:pnb:rrb}
This implementation incorporates the rough relaxed bounding method proposed by Gillard et al. (2021). Rough relaxed bounding was used to trim the domain of each node during construction of the restricted DDs, and was also added as a check to the \emph{filter} function in Algorithm~\ref{algo:RelaxedAlgo}. When the initial model is created, a map is also created from each node $u$, to a list of the other nodes sorted by their distance from $u$. \edit{For the SOP,} the rough relaxed bound $rrb(a)$ of an arc $a_{fg}$ was calculated as follows. For each node $u$ that has not necessarily been visited ($u \notin All_g^\downarrow$), look up the shortest distance from that node to a different node that has also not been visited. Then,  sort the resulting list, and repeatedly remove the largest value until the list has a length equal to the number of remaining decisions. The sum of the values in the list, plus the value of the shortest path from $r$ to the end of $a$, is the rough relaxed bound of $a$. If $rrb(a)$ is worse than the best known solution, the arc is removed. \edit{Slight modifications were made to this procedure to extend it to the TSPTW and Makespan implementations used to generate the results shown in Section \ref{sec:experiments2}}.

\subsection{Peel-and-Bound with Top-Down Compilation}
\label{sec:pnb:topdown}
The structure of peel-and-bound is designed to take advantage of relaxed decision diagrams that are compiled by separation. Here we propose a method for applying peel-and-bound to relaxed decision diagrams that are compiled top-down. However, we have not tested this method, and it remains a topic of future research to determine if it would be useful in practice. The goal of peel-and-bound is to re-use work already done by reusing diagrams. When performing top-down compilation, nodes are merged instead of separated. The peel procedure can be used exactly as before, but after a node is peeled, the remaining diagrams must be altered so that new nodes can be added top-down using a merge procedure. However, there are no nodes to add, the diagrams already represent feasible bounds on the problem, and so some nodes must be removed. Begin by selecting a relaxed node $u$ (in other words some node $u$ that is not exact), and removing it from the diagram. Then remove any arcs in the diagram that are sub-optimal or no longer feasible due to the removed node. For each arc $a_{vu}$ that used to point to $u$, create a new arc $a_{vu\prime}$ that points to a new node $u\prime$ created by following the top-down compilation rules being used. Finally, proceed to perform top-down compilation using the set of new nodes as root nodes for the procedure. 

\section{Initial Experiments on the Sequence Ordering Problem}
\label{sec:results}
The goal of this section is to assess the performances of the peel-and-bound algorithm (PnB, Algorithm~\ref{algo:PnB}) as it compares to the standard decision diagram based branch-and-bound algorithm (BnB, Algorithm~\ref{algo:BnB}). Both algorithms are implemented in Julia and are open-source\footnote{https://github.com/IsaacRudich/PnB$\_$SOP}. 
To ensure a fair comparison, both algorithms resort to the same function for generating relaxed decision diagrams (Algorithm~\ref{algo:RelaxedAlgo}), and the same function for generating restricted decision diagrams. While the functions being called are the same, there are two differences at run-time. At the end of line $26$ in Algorithm~\ref{algo:RelaxedAlgo}, an additional operation runs during BnB where the values of the arcs leaving layer $j$ are updated. The second difference is that BnB starts each relaxation from a $1$-width DD, while PnB passes a partially completed diagram to the relaxation function as a starting point.

The testing environment was built from scratch to ensure a fair comparison, 
so it lacks the many propagators used by cutting-edge solvers like CPO to remove nodes from the PnB/BnB queue~\shortcite{cpS,MDDForSP}. However, it provides a clean comparison of the two algorithms by requiring that every function used by both BnB and PnB is exactly the same between the two, with the only differences arising due to PnB's ability to ensure that all arcs are exact from the beginning. All of the heuristic decisions that were made are identical for both algorithms. These experiments were performed before the implementation of maximal node selection (Section \ref{sec:pnb:implementation:nodeselection}), unlabeled arcs (Section \ref{sec:pnb:implementation:unlabeled}), embedded restricted decision diagrams (Section \ref{sec:pnb:implementation:betterdds}), and the search diversification procedure (Section \ref{sec:pnb:implementation:variety}), so those methods are not included in these experiments. An improved implementation that incorporates those methods and ideas is discussed in Section \ref{sec:experiments2}.

\subsection{Description of the Heuristics Considered}
\label{sec:results:heuristics}
The \emph{sequence ordering problem} can be considered as an asymmetric \emph{travelling salesperson problem} with precedence constraints. The objective is to find a minimum cost path that visits each of the $n$ elements exactly once, and respects the precedence constraints. The method used for generating relaxed DDs requires creating a heuristic ordering of all possible arc assignments by importance. The arc values in this case are representative of the $n$ elements in the path. The ordering used was generated by sorting the $n$ elements, first by their average distance from the other elements, and then by the number of elements each element must precede. The resulting order places a higher importance on elements that are far away from other elements and must precede many other elements.

The branch-and-bound algorithm processes nodes in an order designed to try and improve the existing relaxed bound at each iteration. When a node $u$ is added to the BnB queue, it is assigned a value equal to the value of the shortest path from the root $r$ to the terminal $t$, that passes through $u$. The best known relaxed bound on the problem is the smallest value of a node in the queue, and that node is always chosen to be processed. Peel-and-bound is implemented with the same goal of improving bounds at each iteration. However, PnB stores diagrams, not nodes. Let the value of a diagram be the value of the shortest path to the terminal. At each iteration of peel-and-bound, the diagram with the lowest value is selected, and then a node is chosen from that diagram to induce the peel process. All of the experiments here used a process where the selected node is the first node in the shortest path from $r$ to $t$ with a child node that is not exact (the last exact node). Testing was done to determine whether using the last exact node or the frontier node would perform better for the problem being considered, but there was not a significant difference between the two during any of the tests. Several of the benchmark problems were run using various decision diagram widths, and the last exact node was chosen because it sometimes showed a very slight improvement over the frontier node. While it is likely that this choice makes a difference on some problems, it does not matter for the SOP. 

\subsection{Experimental Results}
\label{sec:results:experiment}
The experiments were performed on a computer equipped with an AMD Rome 7532 at 2.40 GHz with 64Gb RAM. The solver was tested using DD widths of $64, 128$, and $256$ on the $41$ SOP problems available in TSPLIB~\shortcite{tsplib}. For comparisons between PnB and BnB, a timestamp, new bounds, and the length of the remaining queue were recorded each time the bounds on a problem were improved. Another experiment was performed to test the scalability of PnB at width $2048$, for which only the final bounds were recorded. Execution time was limited to $3,600$ seconds.

\begin{figure}[!ht]
    \includegraphics[scale=.53]{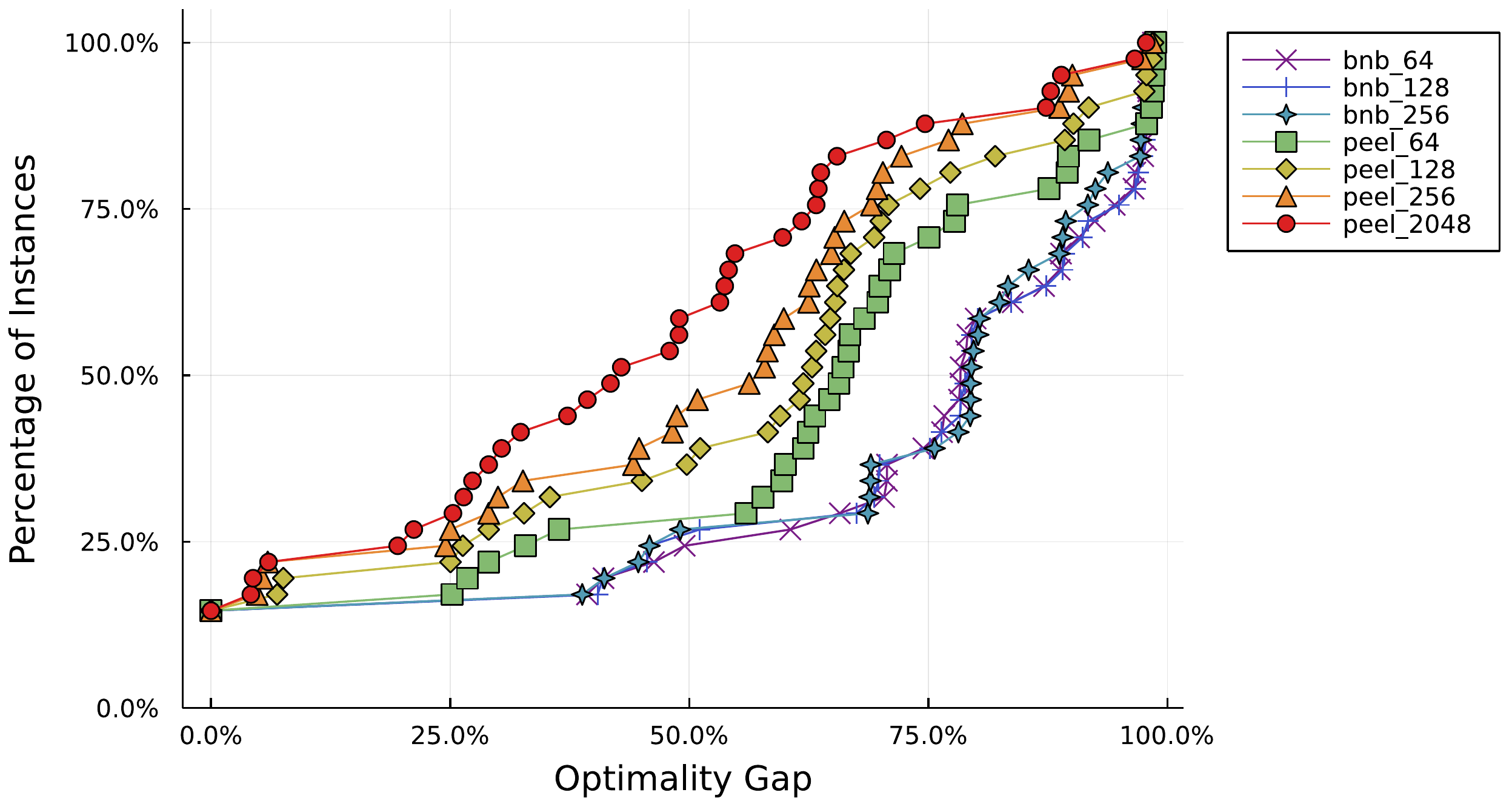}
    \caption{Performance Profiles: the optimality gap $= \frac{upper\_bound - lower\_bound}{upper\_bound}$}
    \label{fig:perfprof2}
\end{figure}

The smallest DD width tested for both methods was $64$, and the largest DD width tested was $256$. \edit{Table} \ref{fig:summary} has summary statistics for those widths as the percentage improvement demonstrated by PnB. A positive percentage always indicates that PnB performed better than BnB in that category, while a negative percentage indicates that BnB performed better. Figure \ref{fig:perfprof2} shows performance profiles for all of the experiments. \edit{Table} \ref{fig:Data4} contains summary statistics comparing PnB at width $256$ to PnB at width $2048$, where a positive percentage always indicates that the width of $2048$ performed better. 

\begin{table}[!ht]
    \centering
    \begin{tabular}{c|c c c c | c c c c|}
        & \multicolumn{4}{c|}{Width: 64} & \multicolumn{4}{c|}{Width: 256} \\
        & RB & BS & OG & QL & RB & BS & OG & QL\\ \hline
        Average \% Improvement & $114\%$ & $0.5\%$ & $22.8\%$ & $1,647\%$ & $545\%$ & $3.3\%$ & $181\%$ & $308\%$\\
        Median \% Improvement & $26\%$ & $0.05\%$ & $17.4\%$ & $734\%$ & $80\%$ & $1.7\%$ & $35\%$ & $141\%$
    \end{tabular}
    \caption{Summary Statistics: percentage improvement of peel-and-bound over branch-and-bound. RB = Relaxed Bound, BS = Best Solution, OG = Optimality Gap, QL = Queue Length. \edit{Tables} \ref{fig:Data1} and \ref{fig:Data2} in Appendix \ref{sec:data} show the comprehensive results.}
    \label{fig:summary}
\end{table}
\begin{table}[!ht]
    \centering
    \begin{tabular}{c|c c c |}
        & \multicolumn{3}{c|}{PnB: 2048 vs PnB: 256}\\
        & Relaxed Bound & Best Solution & Optimality Gap\\ \hline
         Average \% Improvement & $19.5\%$ & $0.8\%$ & $18.6\%$ \\
         Median \% Improvement &  $16.3\%$ & $0.5\%$ & $13.7\%$
    \end{tabular}
    \caption{Summary Statistics: percentage improvement of peel-and-bound at width $2048$ over peel-and-bound at width $256$. Table \ref{fig:Data3} in Appendix \ref{sec:data} shows the comprehensive results.}
    \label{fig:Data4}
\end{table}

As shown in \edit{Table}~\ref{fig:summary}, peel-and-bound vastly outperforms branch-and-bound in these experiments. The average and median improvements from using peel-and-bound at both widths are significant in terms of the relaxed bound, the remaining optimality gap, and the number of nodes that still need to be processed. The best solution found by the end of the runtime also tends to be slightly better with peel-and-bound, but the found solutions are often so close to the real optimal solutions that there is little room for improvement. At both widths, six of the problems were solved to optimality. BnB was faster in only one of those cases, and in that case the difference was $.04$ seconds. The median time for PnB to close in these cases was $191\%$ faster at a width of 64, and $580\%$ faster at a width of $256$. The relaxed bound produced by PnB at a width of $64$ was better for $28$ of the remaining $35$ problems, and at a width of $256$ was better for $34$ of the remaining $35$ problems. The optimality gap was similarly better for peel-and-bound on every problem except the ones where branch-and-bound found a better relaxed bound. However, of the problems where branch-and-bound had a better optimality gap, the improvement was less than $1\%$ for all but one problem. 
Figure~\ref{fig:perfprof2} reinforces that even though there are some instances where a specific branch-and-bound setting slightly outperforms a specific peel-and-bound setting, the gap in those cases is small relative to the general gap between all peel-and-bound settings and all branch-and-bound settings. 

As shown in \edit{Table}~\ref{fig:Data4}, increasing the width to $2048$ from $256$ led to an $19.5\%$ average improvement ($16.3\%$ median improvement) in the relaxed bound. 
Figure~\ref{fig:perfprof2} shows that the performance of peel-and-bound nearly uniformly increases with the maximum allowable width. Similar to the difference between branch-and-bound and peel-and-bound, some specific instances see a small out-performance of the peel-and-bound running at a smaller width, but the gap is small relative to the usual gap between the $2048$-width experiment and the rest of the experiments. 
Additionally, Figure~\ref{fig:perfprof2} shows that peel$\_2048$ solved $50\%$ of instances to within a $42\%$ optimality gap, peel$\_64$ solved $50\%$ of instances to within a $67\%$ optimality gap, and the best performing branch and bound (bnb$\_64$) solved $50\%$ of instances to within only a $79\%$ optimality gap. The overall performance of peel-and-bound improves as more problems are considered, especially as the maximum allowable width for the decision diagrams is increased.

\begin{figure}[!ht]
    \centering
        \setlength{\tabcolsep}{1pt}
        \begin{tabular}{c}
            \includegraphics[scale=.5]{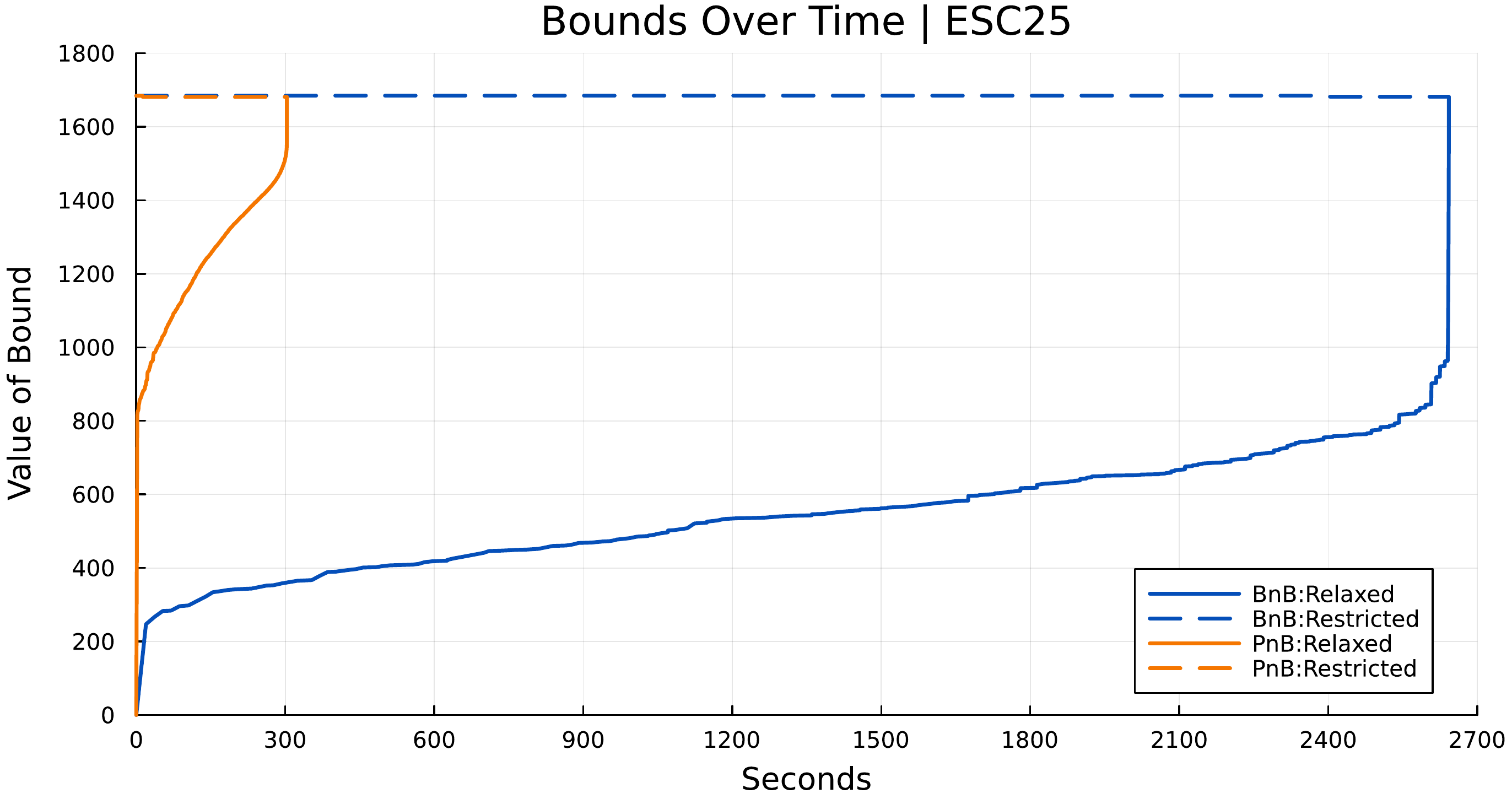} \\
             $\quad$ \\
            \includegraphics[scale=.5]{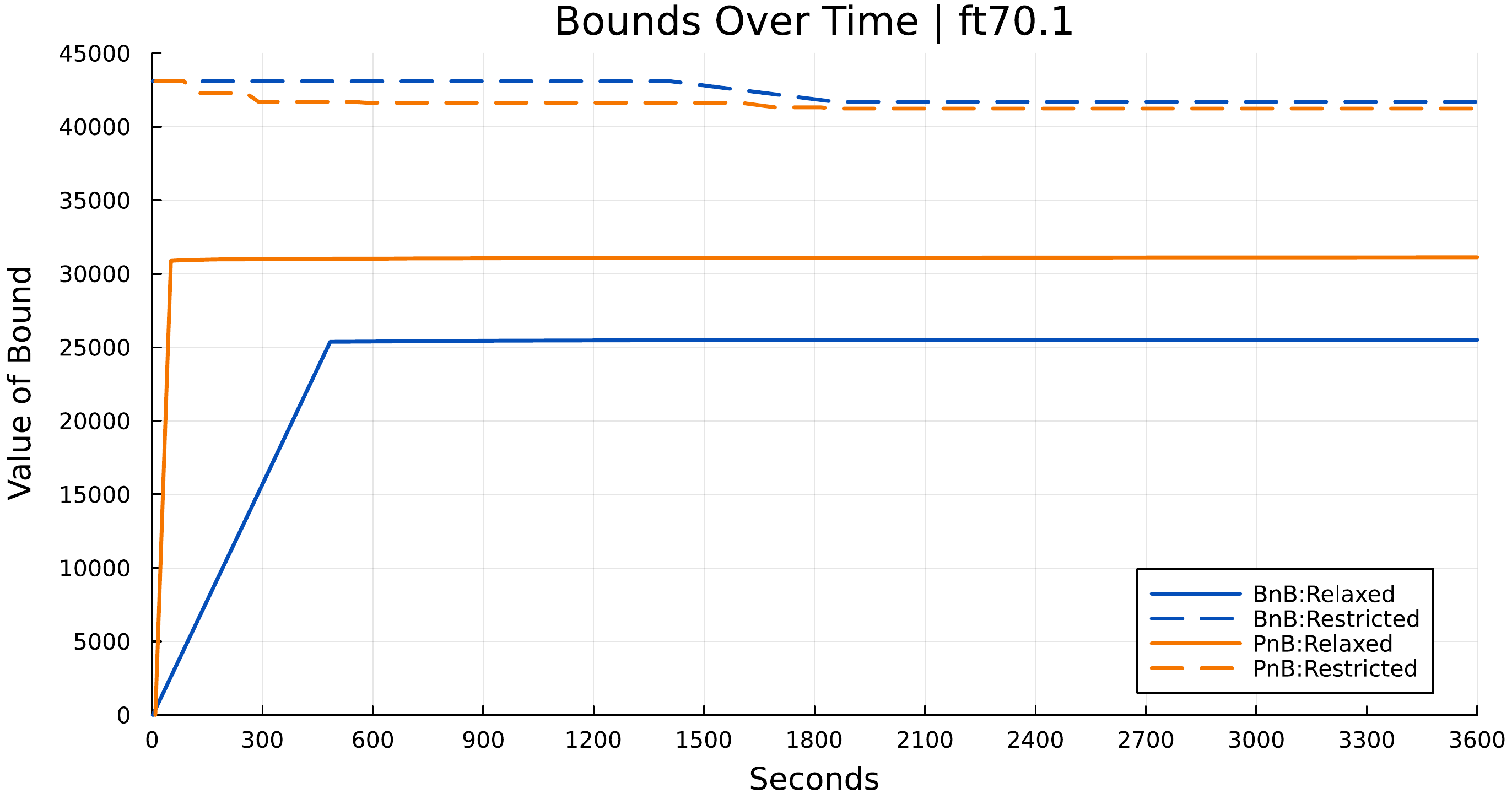}\\
        \end{tabular}
    \caption{\edit{Dual bounds over the runtime of the experiment for ESC25, which was solved within the time limit, and ft70.1, which was not solved within the time limit. Both use a width of 256}.}
    \label{fig:sopExample}
\end{figure}

The selected graphs shown in Figure~\ref{fig:sopExample} are representative of the two main types of behavior observed over the problem set. On problems where the underlying relaxation method works well, the relaxed bound moves quickly towards convergence with the best found solution. On problems where the underlying relaxation does not work well, both algorithms slowly improve the relaxed bound, but PnB starts stronger as it can use exact arc values, and it maintains the advantage throughout. It is clear from the time-series data that to be competitive with cutting-edge solvers, peel-and-bound must be combined with other constraint programming propagators. However, it is also clear that peel-and-bound can have a significant edge over a propagator that generates the required decision diagrams from scratch at each iteration.

\section{Experiments with Improved Methods}
\label{sec:experiments2}
The first implementation of peel-and-bound, which was used to generate the results in Section \ref{sec:results}, was designed to create a fair comparison of peel-and-bound and branch-and-bound; it was also limited to the SOP. We have since re-implemented peel-and-bound so the implementation is generic. The new version is similarly open-source\footnote{https://github.com/IsaacRudich/ImprovedPnB}, and all of the data from the following experiments can be found with the code. In this section, we further explore the performance of the algorithm and compare the performance of a subset of the heuristic methods proposed in Section \ref{sec:pnb}. The experiments were performed on a computer equipped with an AMD Rome 7532 at 2.40 GHz with 186Gb RAM.

\subsection{Node Selection Heuristic}
\label{sec:experiments2:selection}
\begin{figure}[!h]
    \centering
        \setlength{\tabcolsep}{1pt}
        \includegraphics[scale=.5]{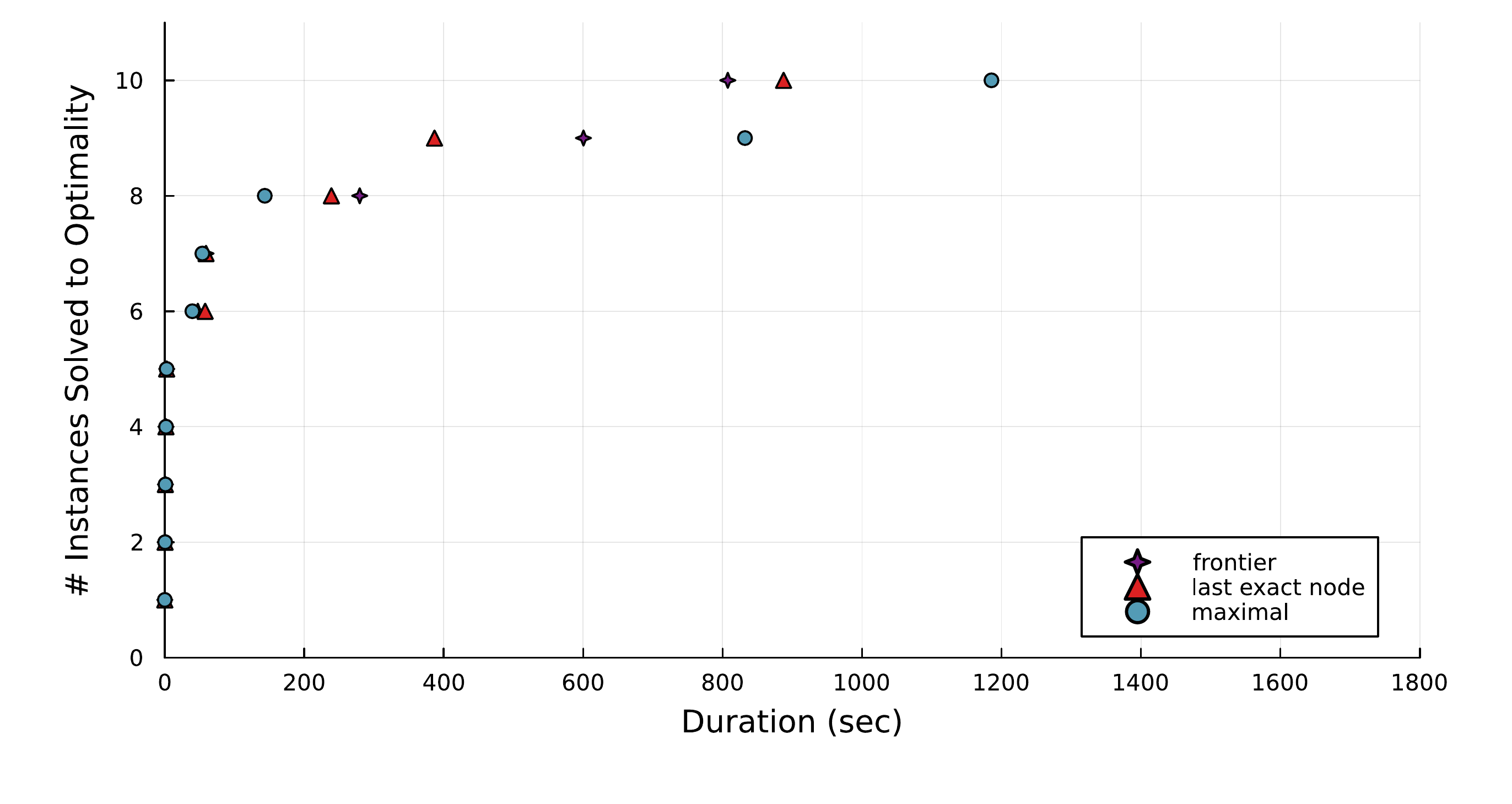}\\
        SOP Instances Solved over Time: max width of 2048\\
        \vspace{20px}
        \includegraphics[scale=.50]{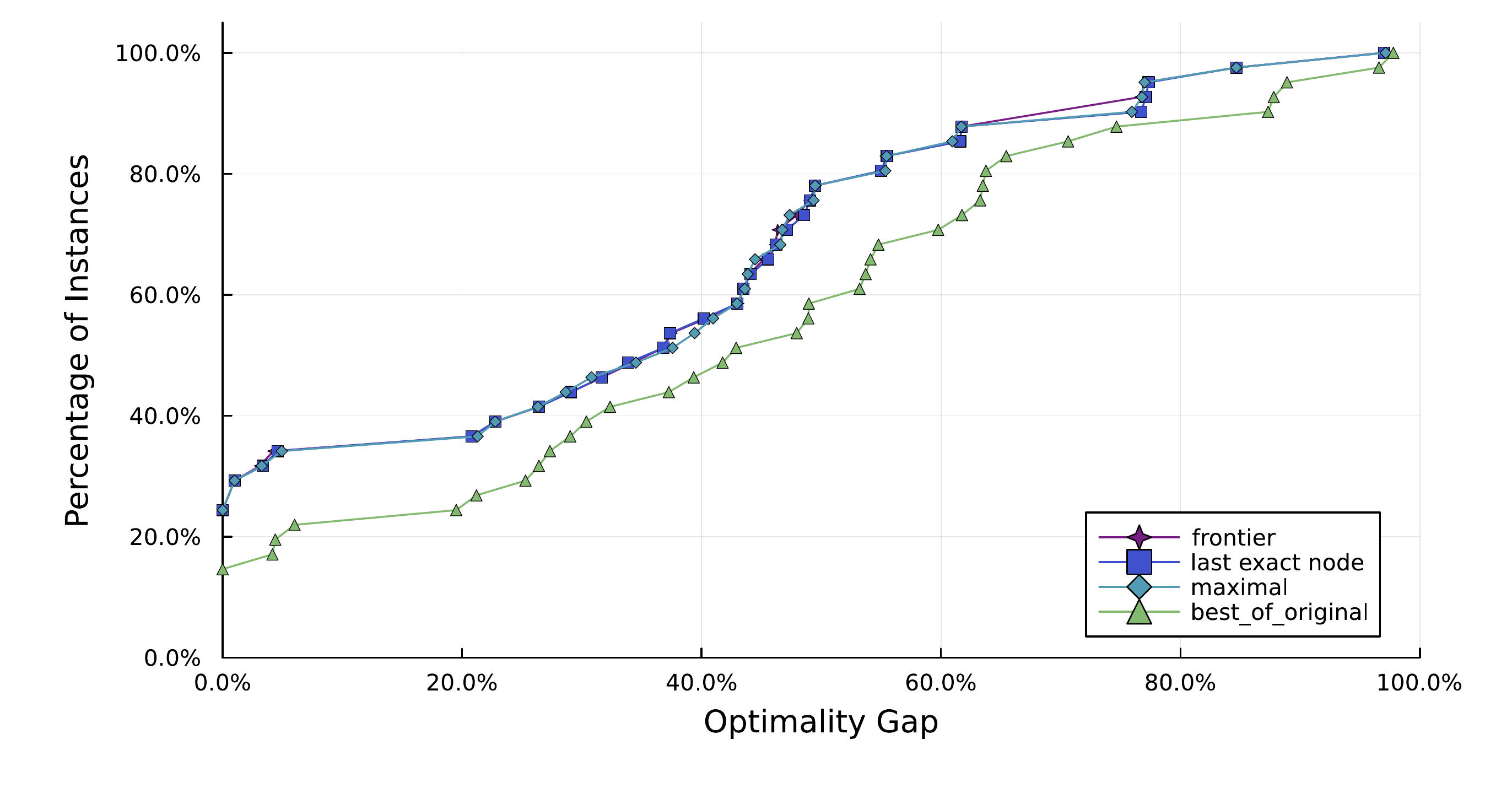}\\
        Performance Profiles on the SOP: max width of 2048\\
    \caption{Performance of Peel-and-Bound on SOP.}
    \label{fig:sopcompare}
\end{figure}

In Section \ref{sec:pnb:implementation:nodeselection}, we propose three heuristics for selecting a node to be peeled: frontier, last exact node, and maximal. Here we compare the performance of those three different settings on the same set of SOP instances we tested in Section \ref{sec:results}; we similarly limit the runtime of the solver to 3600 seconds. We include the results from the most successful run performed by the original implementation of peel-and-bound to show how much the new implementation has improved in general. The results are shown in Figure \ref{fig:sopcompare}. The first graph is a scatter-plot displaying the solve time of each problem that was solved to optimality, for each of the three node selection settings. The second graph shows performance profiles for the same tests, and includes data from the best run of the first implementation of peel-and-bound to show the overall improvement of the solver. All of the new tests were performed using 2048 as the maximum width, and included a diversified search with \edit{$k=5$} using the procedure described in Section \ref{sec:pnb:implementation:variety}.

Figure \ref{fig:sopcompare} makes it clear that choosing between the node selection heuristics has little effect on solving the SOP. The number of instances solved to optimality is identical, with the maximal having a slight lag on the final two problems, and the performance profiles are nearly identical. The performance profiles also serve to \edit{demonstrate} the progress of the solver, with the percentage of instances closed to any given optimality gap being about 10\% higher. 

\subsection{Traveling Salesman Problem with Time Windows}
\label{sec:experiments2:tsptw}
\begin{figure}[!h]
    \centering
        \setlength{\tabcolsep}{1pt}
        \includegraphics[scale=.5]{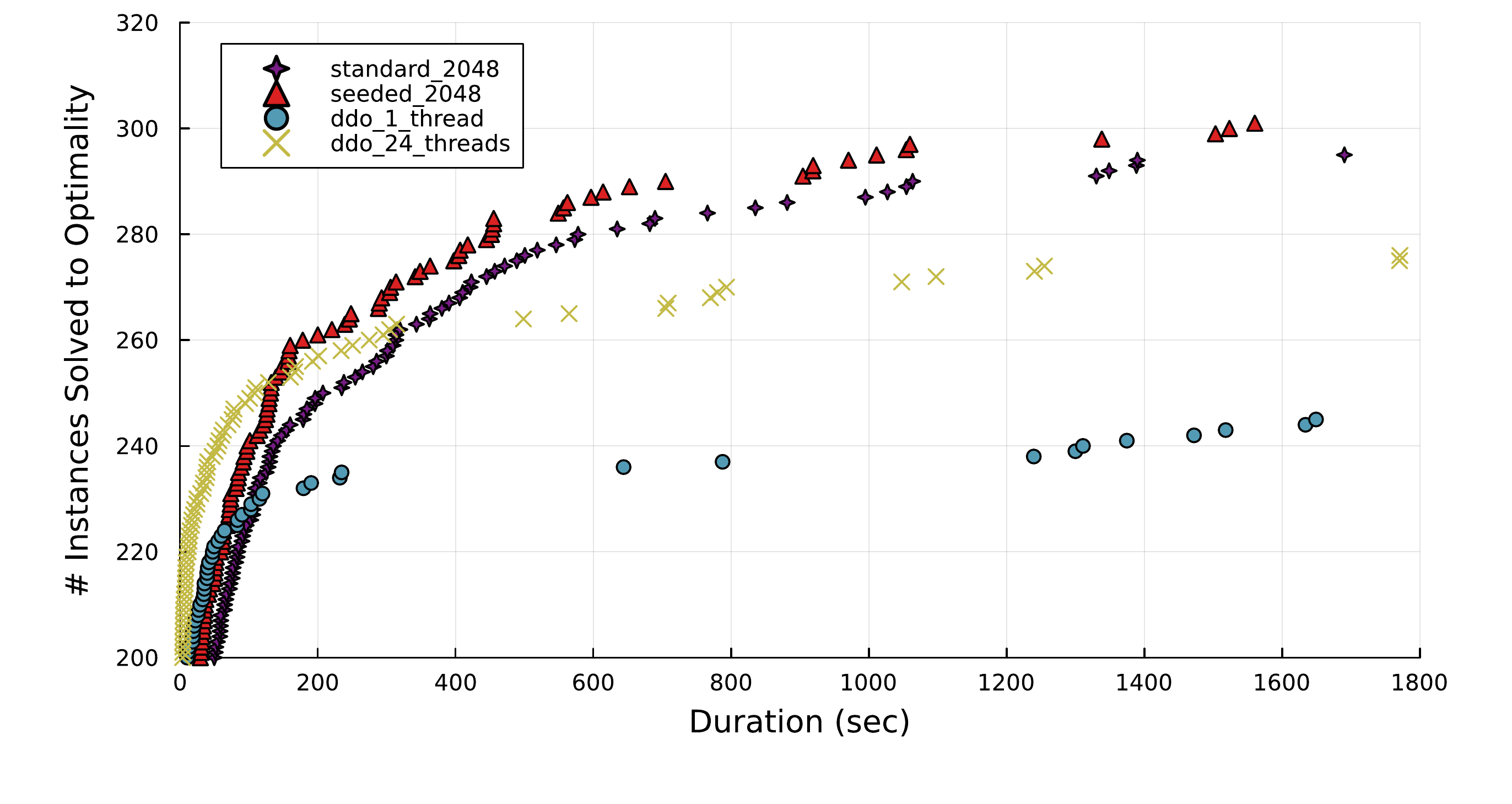}\\
        TSPTW Instances Solved over Time: DDO vs. Peel-and-Bound\\
        \vspace{20px}
        \includegraphics[scale=.5]{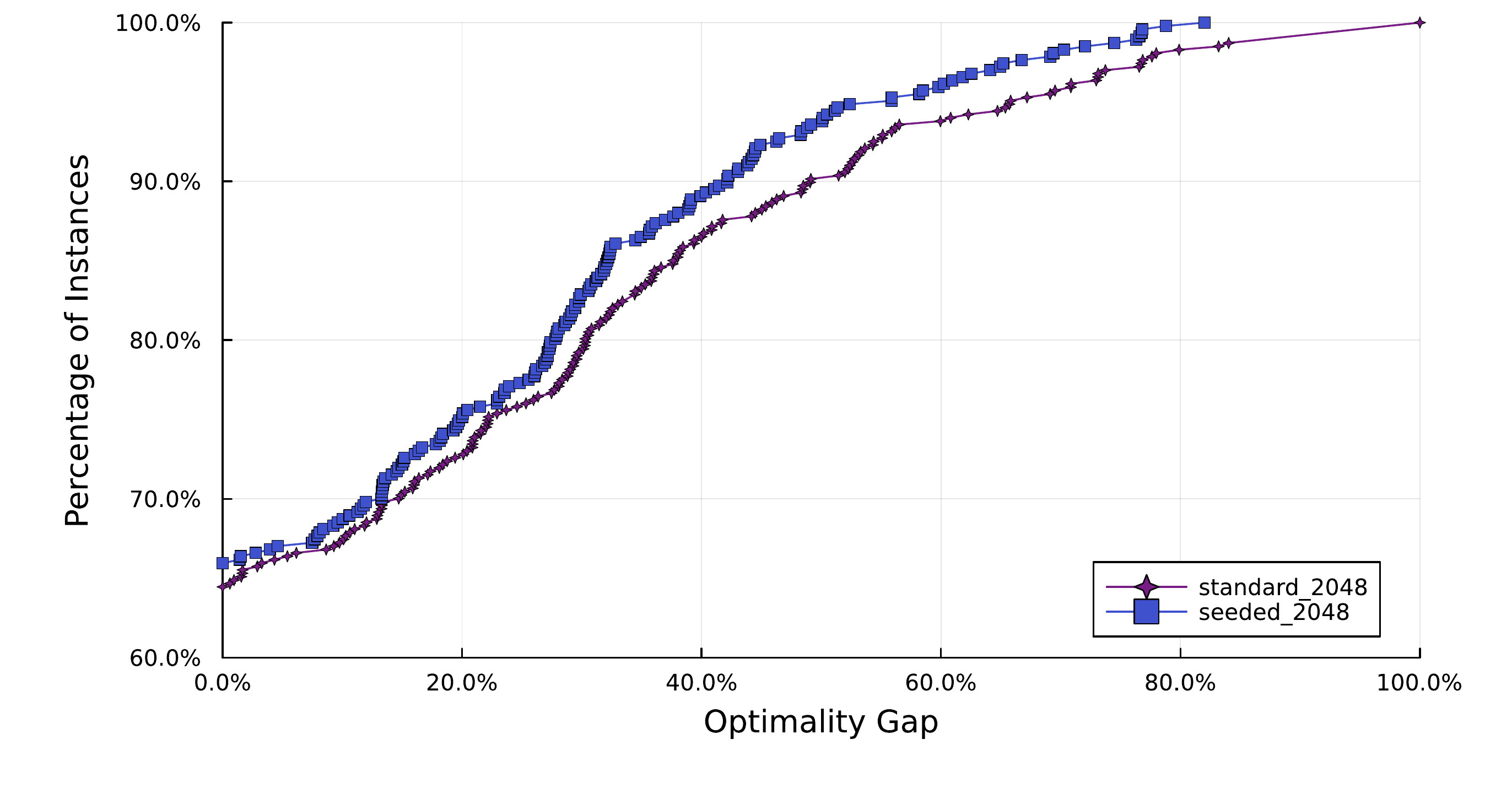}\\
        TSPTW Performance Profiles\\
    \caption{Performance of Peel-and-Bound on TSPTW.}
    \label{fig:tsptw}
\end{figure}

The traveling salesman problem with time windows (TSPTW) is a variation of the traveling salesman problem where a salesman must find the shortest cycle that visits each of his customers once, but each customer may have an added constraint requiring they be visited within a specific time window. It provides a useful metric for benchmarking the performance of the implementation of peel-and-bound, as decision diagram focused algorithms are more likely to outperform traditional methods on highly constrained problems than highly unconstrained problems, and the TSPTW instances are highly constrained. For the following tests we include both a standard and seeded run of the solver, where seeded means that the solver started out knowing the value of the best known solution, and skipped the initial diversified search. \edit{The seeded version is of interest because the solver can leverage heuristically generated solutions to reach proof of optimally faster. The difference between the standard and seeded run help may help to determine if taking that step is worthwhile when solving a specific problem.} We test our solver on the same set of $467$ benchmark instances used by Gillard et al. (2021) to test their implementation of decision diagram based branch-and-bound (ddo). These instances are available from L\'opez-Ib\'a\~nez \& Blum (2022), and include the following sets: AFG \shortcite{AFG}, Dumas \shortcite{Dumas}, GendreauDumas \shortcite{gendreaudumas}, Langevin \shortcite{langevin}, Ohlmann-Thomas \shortcite{ohlmann}, Solomon-Pesant \shortcite{pesant1998exact}, and Solomon-Potvin-Bengio \shortcite{potvin1996vehicle}. Figure \ref{fig:tsptw} shows a time to solve graph for the closed instances, once with peel-and-bound on a single thread, once with seeded peel-and-bound on a single thread, once with ddo on a single thread, and once with ddo using 24 threads. The experiments were limited to 3600 seconds. Figure \ref{fig:tsptw} also shows the performance profiles of peel-and-bound using the same data. Since the solution space of TSPTW tends to be drastically more constrained than SOP, at least for the benchmark instances, it is more difficult to find feasible solutions, and we use a width of 100 for the initial diversified search.

It is clear from Figure \ref{fig:tsptw}, that peel-and-bound outperforms ddo on the TSPTW benchmark set. Although ddo is faster to start, peel-and-bound on a single thread solves about 60 more instances than ddo on a single thread, and about 20 more instances than ddo using 24 threads. This experiment plainly demonstrates the advantages of trading memory for speed when using decision diagrams. The performance profiles show that peel-and-bound solved about $65\%$ of the instances to optimality. The seeded version of the solver performed only slightly better. Looking at the raw data (available in the repository), it is also clear that performance degrades when peel-and-bound has trouble finding good solutions, and when the size of the time windows is large (causing the problem to be more unconstrained than others in the benchmark set).

To the best of our knowledge, the last paper to report relaxed bounds on the instances in these benchmark sets is Baldacci et al. (2012). So we compare our results to Baldacci et al. (2012). The standard run of the solver closes \edit{14} instances that were left open by Baldacci et al. (2012), and the seeded run of the solver closes one additional instance. The results for these instances, and the best relaxed bound found by Baldacci et al. (2012), are reported in \edit{Table} \ref{fig:tsptwdata}. A table with every problem from the benchmark set that remains open to the best of our knowledge (in other words, those unsolved by both Baldacci et al. (2012) and peel-and-bound), is available in Appendix \ref{sec:data}: \edit{Table} \ref{fig:tsptwopendata}. The full data is available in Appendix \ref{sec:data}: \edit{Tables} \ref{fig:tsptwcloseddata}-\ref{fig:tsptwcloseddata6}. The raw data is also available in the repository with the solver.

\begin{table}[!h]
    \centering
    \resizebox{\textwidth}{!}{%
        \begin{tabular}{|ccc|c|cccc|cccc|}
            \multicolumn{3}{|c|}{Problem Information} & \multicolumn{1}{l|}{Baldacci et al. (2012)} & \multicolumn{4}{c|}{Single Thread: 2048} & \multicolumn{4}{c|}{Seeded   Single Thread: 2048} \\
            Set & Name & Best Known Solution & LB & LB & UB & Time & OG & LB & UB & Time & OG \\ \hline
            \multirow{2}{*}{AFG}
             & rbg086a.tw & 8400 & 8399 & 8400 & 8400 & 254.47 & - & 8400 & 8400 & 132.49 & - \\
             & rbg092a.tw & 7158 & 7156.6 & 7158 & 7158 & 1917.74 & - & 7158 & 7158 & 288.02 & - \\ \hline
            \multirow{3}{*}{GendreauDumas} & n20w200.004.txt & 293 & 289.484 & 293 & 293 & 27.97 & - & 293 & 293 & 15.47 & - \\
             & n40w200.002.txt & 303 & 302.091 & 303 & 303 & 1388.47 & - & 303 & 303 & 417.75 & - \\
             & n80w100.004.txt & 649 & 645.408 & 649 & 649 & 2372.02 & - & 649 & 649 & 904.33 & - \\ \hline
            SolomonPesant & rc203.1 & 726.99 & 726.66 & 726.99 & 726.99 & 299.58 & - & 726.99 & 726.99 & 200.01 & - \\ \hline
            \multirow{9}{*}{SolomonPotvinBengio} & rc\_202.3.txt & 837.72 & 835.87 & 837.72 & 837.72 & 9.66 & - & 837.72 & 837.72 & 8.67 & - \\
             & rc\_202.4.txt & 793.03 & 791.54 & 793.03 & 793.03 & 108.08 & - & 793.03 & 793.03 & 73.79 & - \\
             & rc\_205.4.txt & 760.47 & 756.95 & 760.47 & 760.47 & 13.65 & - & 760.47 & 760.47 & 9.74 & - \\
             & rc\_206.2.txt & 828.06 & 826.66 & 828.06 & 828.06 & 102.63 & - & 828.06 & 828.06 & 68.91 & - \\
             & rc\_206.4.txt & 831.67 & 827.54 & 831.67 & 831.67 & 107.39 & - & 831.67 & 831.67 & 84.51 & - \\
             & rc\_207.1.txt & 732.68 & 731.57 & 732.68 & 732.68 & 195.9 & - & 732.68 & 732.68 & 157.17 & - \\
             & rc\_207.2.txt & 701.25 & 694.22 & 701.25 & 701.25 & 1690.43 & - & 701.25 & 701.25 & 92.84 & - \\
             & rc\_207.3.txt & 682.40 & 677.23 & 682.40 & 682.40 & 1054.62 & - & 682.40 & 682.40 & 596.7 & - \\
             & rc\_208.1.txt & 789.25 & 785.69 & 751.20 & 794.17 & - & 5.41 & 789.25 & 789.25 & 2511.58 & - \\ \hline
            \end{tabular}
    }
    \caption{TSPTW Results for Newly Closed Problems: LB = Lower Bound, UB = Upper Bound, OG = Optimality Gap. \edit{In all cases the existing best known solution was proven optimal by Peel and Bound.}}
    \label{fig:tsptwdata}
\end{table}

\subsection{Traveling Salesman Problem with Time Windows - Makespan}
\label{sec:experiments2:makepspan}
\begin{figure}[!h]
    \centering
        \setlength{\tabcolsep}{1pt}
        \includegraphics[scale=.5]{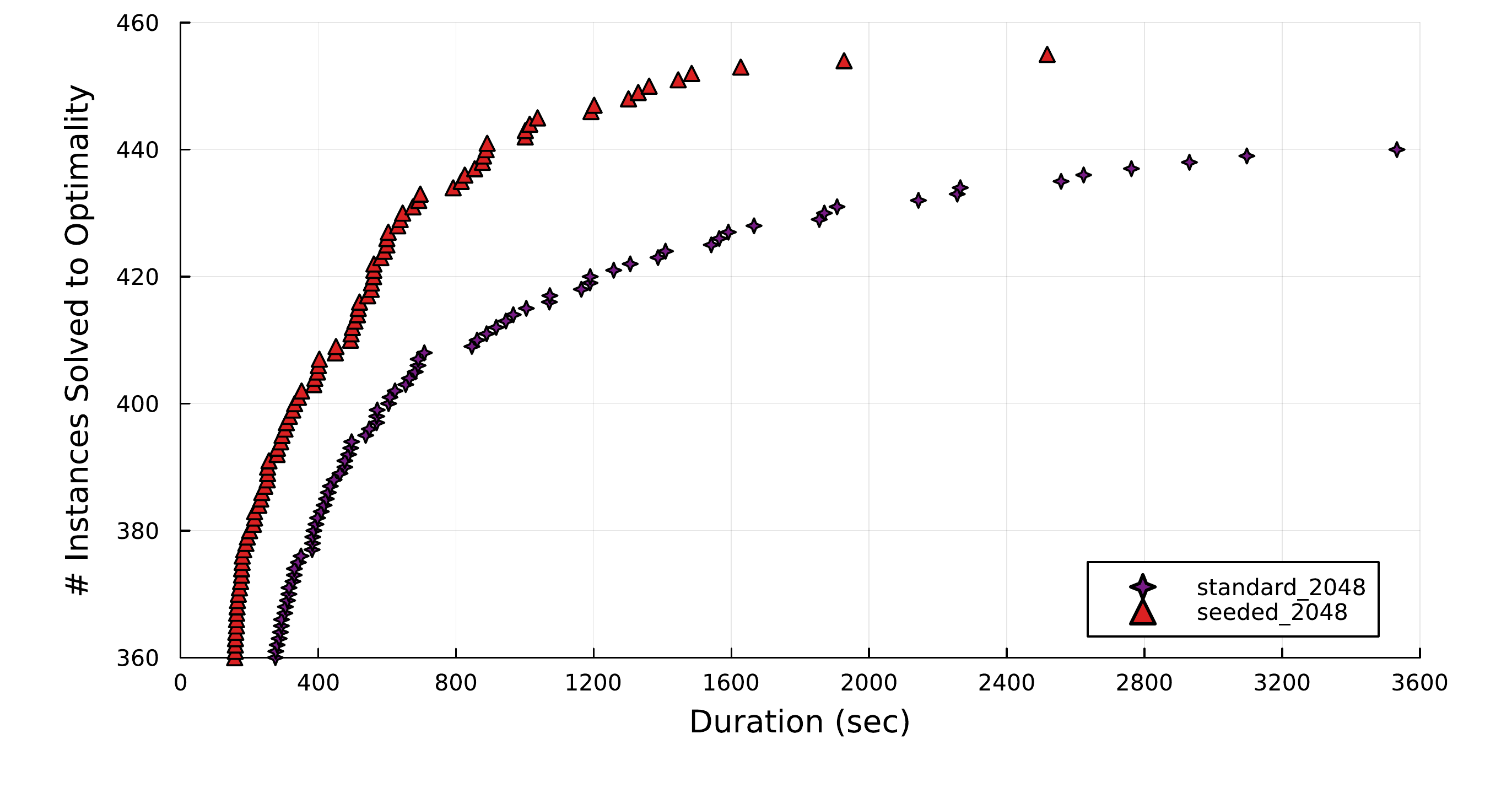}\\
        Makespan Instances Solved over Time\\
        \vspace{20px}
        \includegraphics[scale=.5]{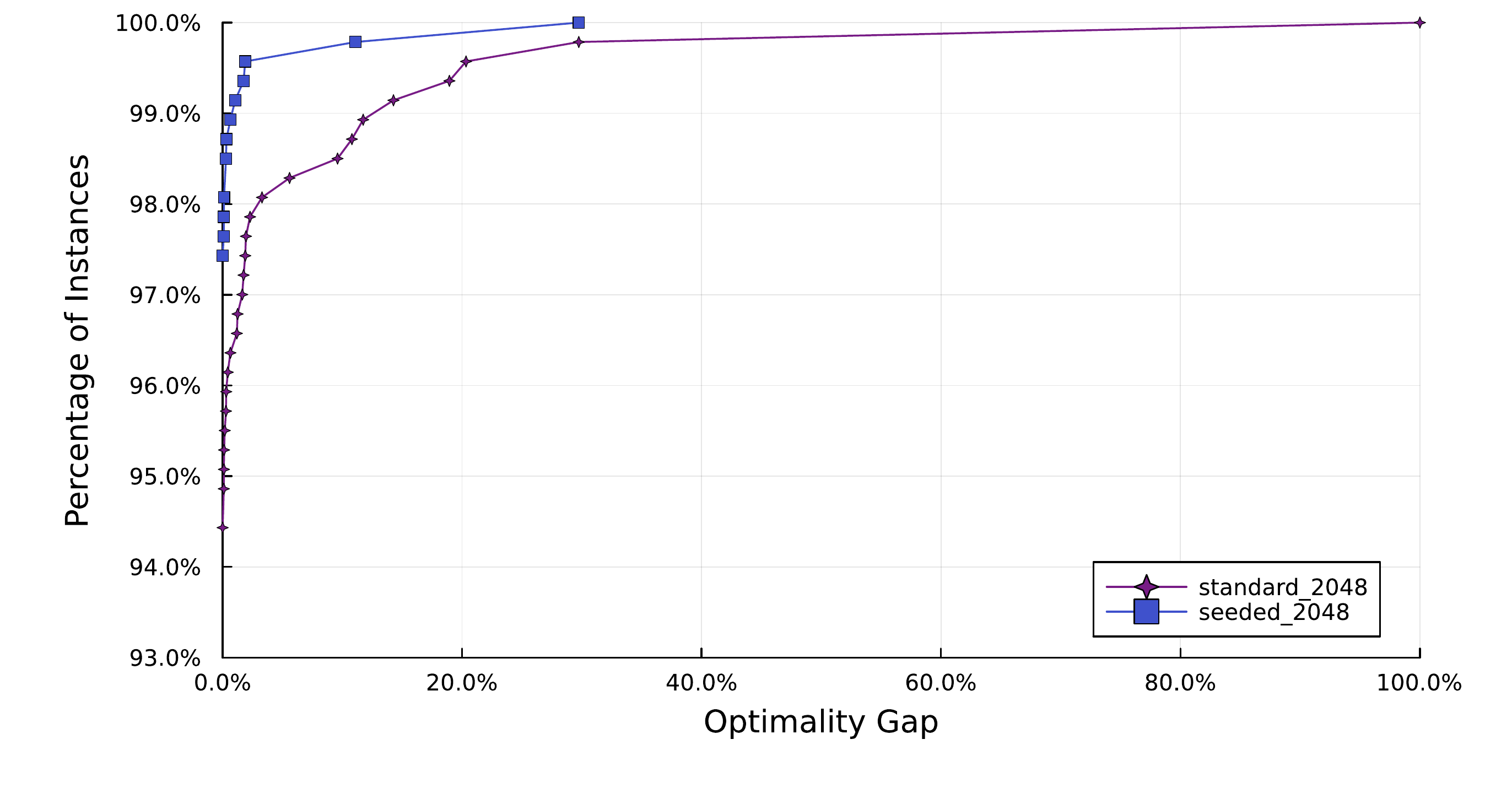}\\
        Makespan Performance Profiles\\
    \caption{Performance of Peel-and-Bound on Makespan.}
    \label{fig:makespan}
\end{figure}

\edit{Makespan adjusts the objective function of TSPTW so} the total time includes any idle time spent waiting for a customer to be available when the salesman arrives early, as opposed to just the distance traveled. To test this problem, we use the same benchmark instances that we did in Section \ref{sec:experiments2:tsptw}, and simply adjust the objective function. The tests all use the same settings as the test for TSPTW, and the results are presented in the same way in Figure \ref{fig:makespan}.

While we were unable to find a a benchmark comparison in the literature with relaxed bounds for the makespan variant of these benchmark sets, the results from testing peel-and-bound on makespan speak for themselves. Just over $94\%$ of the instances were closed to optimality by the unseeded run of the solver, and just over $97\%$ were closed by the seeded run. For the unseeded run, just over $97.5\%$ were closed to $1\%$ optimality gap, and $99\%$ for the seeded run. In total, 26 instances were not closed by the unseeded run, and 12 of those were not closed by the seeded run. The results and bounds for those 26 instances are shown in \edit{Table} \ref{fig:makespandata}. We are unsure if any of the makespan problems are considered to be open, but we can say that the only ones that aren't definitively closed are the 12 \edit{not in bold} that were not solved during the seeded run of the solver. Of those, only 2 have an optimality gap larger than $2\%$. For the 455 closed instances, the best solution found was the best solution reported by L\'opez-Ib\'a\~nez \& Blum (2022). A table with full results for the benchmark sets, and time to solve for all of the closed problems, is shown in Appendix \ref{sec:data}: \edit{Tables} \ref{fig:makespanfulldata}-\ref{fig:makespanfulldata7}. The raw data is also available in the repository with the solver.

\begin{table}[!h]
    \centering
    \resizebox{\textwidth}{!}{%
        \begin{tabular}{|ccc|ccc|ccc|}
        \multicolumn{3}{|c|}{Problem Information} & \multicolumn{3}{c|}{Single Thread: 2048} & \multicolumn{3}{c|}{Seeded Single Thread: 2048} \\
        Set & Name & Best Known Solution & LB & UB & OG & LB & UB & OG \\ \hline
        \multirow{2}{*}{AFG} & rbg050b.tw & 11957 & 11748 & 11957 & 1.75 & 11748 & 11957 & 1.75 \\
         & rbg172a.tw & 17783 & 17766 & 17784 & 0.10 & 17766 & 17783 & 0.10 \\ \hline
        \multirow{2}{*}{Dumas} & \textbf{n150w60.002.txt} & \textbf{940} & \textbf{940} & \textbf{941} & \textbf{0.11} & \textbf{940} & \textbf{940} & \textbf{-} \\
         & \textbf{n200w40.002.txt} & \textbf{1137} & \textbf{1137} & \textbf{-} & \textbf{100} & \textbf{1137} & \textbf{1137} & \textbf{-} \\ \hline
        \multirow{8}{*}{GendreauDumas} & n100w100.003.txt & 819 & 818 & 819 & 0.12 & 818 & 819 & 0.12 \\
         & \textbf{n100w80.003.txt} & \textbf{829} & \textbf{829} & \textbf{839} & \textbf{1.19} & \textbf{829} & \textbf{829} & \textbf{-} \\
         & \textbf{n60w200.003.txt} & \textbf{560} & \textbf{560} & \textbf{561} & \textbf{0.18} & \textbf{560} & \textbf{560} & \textbf{-} \\
         & n80w140.005.txt & 739 & 725 & 739 & 1.89 & 725 & 739 & 1.89 \\
         & \textbf{n80w160.002.txt} & \textbf{654} & \textbf{652} & \textbf{665} & \textbf{1.95} & \textbf{654} & \textbf{654} & \textbf{-} \\
         & n80w180.002.txt & 633 & 631 & 639 & 1.25 & 631 & 633 & 0.32 \\
         & n80w180.005.txt & 632 & 444 & 632 & 29.75 & 444 & 632 & 29.75 \\
         & n80w200.004.txt & 667 & 660 & 671 & 1.64 & 660 & 667 & 1.05 \\ \hline
        \multirow{11}{*}{OhlmannThomas} & \textbf{n150w120.004.txt} & \textbf{925} & \textbf{925} & \textbf{929} & \textbf{0.43} & \textbf{925} & \textbf{925} & \textbf{-} \\
         & \textbf{n150w140.002.txt} & \textbf{1020} & \textbf{1020} & \textbf{1021} & \textbf{0.10} & \textbf{1020} & \textbf{1020} & \textbf{-} \\
         & \textbf{n150w140.004.txt} & \textbf{898} & \textbf{898} & \textbf{919} & \textbf{2.29} & \textbf{898} & \textbf{898} & \textbf{-} \\
         & \textbf{n150w140.005.txt} & \textbf{926} & \textbf{826} & \textbf{926} & \textbf{10.8} & \textbf{926} & \textbf{926} & \textbf{-} \\
         & \textbf{n150w160.002.txt} & \textbf{890} & \textbf{861} & \textbf{912} & \textbf{5.59} & \textbf{890} & \textbf{890} & \textbf{-} \\
         & \textbf{n150w160.004.txt} & \textbf{912} & \textbf{912} & \textbf{943} & \textbf{3.29} & \textbf{912} & \textbf{912} & \textbf{-} \\
         & n200w120.001.txt & 1089 & 1086 & 1089 & 0.28 & 1086 & 1089 & 0.28 \\
         & n200w120.002.txt & 1072 & 1065 & 1072 & 0.65 & 1065 & 1072 & 0.65 \\
         & \textbf{n200w140.001.txt} & \textbf{1138} & \textbf{929} & \textbf{1146} & \textbf{18.94} & \textbf{1138} & \textbf{1138} & \textbf{-} \\
         & n200w140.003.txt & 1083 & 979 & 1083 & 9.60 & 1082 & 1083 & 0.09 \\
         & \textbf{n200w140.005.txt} & \textbf{1121} & \textbf{961} & \textbf{1121} & \textbf{14.27} & \textbf{1121} & \textbf{1121} & \textbf{-} \\ \hline
        SolomonPesant & \textbf{rc204.2} & \textbf{870.52} & \textbf{728.94} & \textbf{914.89} & \textbf{20.33} & \textbf{870.52} & \textbf{870.52} & \textbf{-} \\ \hline
        \multirow{2}{*}{SolomonPotvinBengio} & rc\_204.1.txt & 917.83 & 915.22 & 918.01 & 0.30 & 915.22 & 917.83 & 0.28 \\
         & rc\_208.3.txt & 686.80 & 606.15 & 686.80 & 11.74 & 610.58 & 686.80 & 11.10
        \end{tabular}
    }
    \caption{Makespan Results for \edit{Problems Not Closed by the Unseeded Run: LB = Lower Bound, UB = Upper Bound, OG = Optimality Gap. The bolded instances were closed by the seeded run.}}
    \label{fig:makespandata}
\end{table}
\newpage
\section{Conclusion and Future Work}
\label{sec:conclusion}
This paper presented a peel-and-bound algorithm as an alternative to branch-and-bound. In peel-and-bound, constructed decision diagrams are repeatedly reused to avoid unnecessary computation. Additionally, peel-and-bound can be used in combination with a decision diagram structure that only admits exact arc values, to increase scalability and further reduce the amount of work needed at each iteration of the algorithm. We identified several heuristic decisions that can be used to adjust peel-and-bound, and provided insight into how the algorithm can be applied to other problems.

We compared the performance of a peel-and-bound scheme to a branch-and-bound scheme using the same DD based propagator. We tested both algorithms on the 41 instances of the SOP from TSPLIB. Results show that peel-and-bound significantly outperforms branch-and-bound on the SOP by generating substantially stronger relaxed bounds on instances that were not closed during the experiment, and reaching optimality faster when the instances were closed. We then re-implemented the algorithm to be more efficient, and generic. We tested the new implementation on the 467 benchmark instances of TSPTW used by Gillard et al. (2021) to test their decision diagram based branch-and-bound solver (ddo). The results show that peel-and-bound outperforms ddo on TSPTW, even when ddo is using parallel processing. Furthermore, peel-and-bound closed 16 instances that, to the best of our knowledge, are open in the literature. In our final test, we ran the new implementation of peel-and-bound on the makespan variant of the 467 TSPTW instances. Peel-and-bound closed $94\%$ of the makespan instances, and an additional $3\%$ when seeded with the best known solution. We provide best known bounds for all TSPTW and makespan instances that we believe to be open. 

This paper provides strong support for the value of re-using work in DD based solvers. Furthermore, peel-and-bound benefits from scaling the maximum allowable width. Thus, relaxed DDs that yield strong bounds at scale, but are too costly to generate iteratively, only need to be constructed once. The method detailed in this paper focused on DDs generated by separation; future research could focus on DDs generated using a merge operator.

\appendix

\section{Experimental Data}
\label{sec:data}

\begin{table}[!h]
    \centering
    \resizebox{!}{5.96cm}{%
\begin{tabular}{|cc|ccccc|ccccc|ccccc|}
        \multicolumn{2}{c|}{Problem Info} & \multicolumn{5}{c|}{BnB: width 64} & \multicolumn{5}{c|}{PnB: width 64} & \multicolumn{5}{c}{Percent Improvements} \\
        Name & n & \multicolumn{1}{c}{RB} & \multicolumn{1}{c}{BS} & T & OG & QL & \multicolumn{1}{c}{RB} & \multicolumn{1}{c}{BS} & T & OG & QL & \multicolumn{1}{c}{RB} & BS & T & OG & QL \\ \hline
        ESC07 & 9 & 2,125 & 2,125 & 0.03 & 0\% & - & 2,125 & 2,125 & 0.07 & 0\% & - &  &  & -57\% &  &  \\
        ESC11 & 13 & 2,075 & 2,075 & 0.65 & 0\% & - & 2,075 & 2,075 & 0.42 & 0\% & - &  &  & 55\% &  &  \\
        ESC12 & 14 & 1,675 & 1,675 & 1.99 & 0\% & - & 1,675 & 1,675 & 0.64 & 0\% & - &  &  & 211\% &  &  \\
        ESC25 & 27 & 1,681 & 1,681 & 956 & 0\% & - & 1,681 & 1,681 & 353 & 0\% & - &  &  & 171\% &  &  \\
        ESC47 & 49 & 334 & 1,542 &  & 78\% & 8,842 & 368 & 1,676 &  & 78\% & 1,295 & 10.2\% & -8.0\% &  & 0.4\% & 583\% \\
        ESC63 & 65 & 8 & 62 &  & 87\% & 2,756 & 44 & 62 &  & 29\% & 15 & 450.0\% & 0.0\% &  & 200.0\% & 18273\% \\
\end{tabular}
    }
     \caption{Comparison Data for width 64 experiments on SOP: RB = Relaxed Bound, BS = Best Solution, T = Time in Seconds, OG = Optimality Gap, QL = Queue Length. Full time series data is available in the GitHub repository.}
     \label{fig:Data1}
\end{table}

\begin{table}[!ht]
    \centering
    \resizebox{!}{5.96cm}{%
    \begin{tabular}{|cc|ccccc|ccccc|ccccc}
        \multicolumn{2}{c|}{Problem   Info} & \multicolumn{5}{c|}{BnB: width 256} & \multicolumn{5}{c|}{PnB: width 256} & \multicolumn{5}{c}{Percent Improvements} \\
        Name & n & \multicolumn{1}{c}{RB} & \multicolumn{1}{c}{BS} & T & OG & QL & \multicolumn{1}{c}{RB} & \multicolumn{1}{c}{BS} & T & OG & QL & \multicolumn{1}{c}{RB} & BS & T & OG & QL \\ \hline
        ESC07 & 9 & 2,125 & 2,125 & 0.04 & 0\% & - & 2,125 & 2,125 & 0.04 & 0\% & - &  &  & 0\% &  &  \\
        ESC11 & 13 & 2,075 & 2,075 & 0.48 & 0\% & - & 2,075 & 2,075 & 0.41 & 0\% & - &  &  & 17\% &  &  \\
        ESC12 & 14 & 1,675 & 1,675 & 1.66 & 0\% & - & 1,675 & 1,675 & 0.34 & 0\% & - &  &  & 388\% &  &  \\
        ESC25 & 27 & 1,681 & 1,681 & 2,643 & 0\% & - & 1,681 & 1,681 & 303 & 0\% & - &  &  & 771\% &  &  \\
        ESC47 & 49 & 312 & 1,590 &  & 80\% & 720 & 658 & 1,339 &  & 51\% & 740 & 110.9\% & 18.7\% &  & 58.0\% & -3\% \\
        ESC63 & 65 & 9 & 62 &  & 85\% & 53 & 44 & 62 &  & 29\% & 3 & 388.9\% & 0.0\% &  & 194.4\% & 1667\% \\
\end{tabular}
    }
     \caption{Comparison data for width 256 experiments on SOP: RB = Relaxed Bound, BS = Best Solution, T = Time in Seconds, OG = Optimality Gap, QL = Queue Length. Full time series data is available in the GitHub repository.}
     \label{fig:Data2}
 \end{table}

\begin{table}[!ht]
    \centering
    \resizebox{!}{5.96cm}{%
    \begin{tabular}{|cc|ccc|ccc|ccc|}
        \multicolumn{2}{|c|}{Problem   Info} & \multicolumn{3}{c|}{PnB: width 256} & \multicolumn{3}{c|}{PnB: width 2048} & \multicolumn{3}{c|}{Percent Improvements} \\
        Name & n & \multicolumn{1}{c}{RB} & \multicolumn{1}{c}{BS} & OG & \multicolumn{1}{c}{RB} & \multicolumn{1}{c}{BS} & OG & \multicolumn{1}{c}{RB} & BS & OG \\ \hline
        ESC47 & 49 & 658 & 1,339 & 51\% & 882 & 1,304 & 32\% & 34.0\% & 2.7\% & 57.2\% \\
        ESC63 & 65 & 44 & 62 & 29\% & 44 & 62 & 29\% & 0.0\% & 0.0\% & 0\% \\
        ESC78 & 80 & 5,600 & 20,135 & 72\% & 6,025 & 20,505 & 71\% & 7.6\% & -1.8\% & 2.2\% \\
        ft53.1 & 54 & 4,603 & 8,244 & 44\% & 5,167 & 8,237 & 37\% & 12.3\% & 0.1\% & 18.5\% \\
        ft53.2 & 54 & 3,555 & 8,648 & 59\% & 4,910 & 8,598 & 43\% & 38.1\% & 0.6\% & 37.3\% \\
        ft53.3 & 54 & 4,852 & 11,095 & 56\% & 7,722 & 11,092 & 30\% & 59.2\% & 0.0\% & 85.2\% \\
\end{tabular}
    }
    \caption{Comparison of PnB at 2048 over PnB at 256 on SOP: RB = Relaxed Bound, BS = Best Solution, OG = Optimality Gap.}
    \label{fig:Data3}
 \end{table}

 \begin{table}[!ht]
    \centering
    \resizebox{!}{6.79cm}{%
            \begin{tabular}{|ccc|c|ccc|ccc|}
                \multicolumn{3}{|c|}{Problem Information} & \multicolumn{1}{l|}{Baldacci et. al.   (2012)} & \multicolumn{3}{c|}{Single Thread: 2048} & \multicolumn{3}{c|}{Seeded Single Thread:   2048} \\
                Set & Name & Best Known Solution & LB & LB & UB & OG & LB & UB & OG \\ \hline
                Dumas & n100w60.001.txt & 655 & NA & 602 & 659 & 8.65 & 602 & 655 & 8.09 \\
                Dumas & n100w60.003.txt & 744 & NA & 619 & 749 & 17.36 & 620 & 744 & 16.67 \\
                Dumas & n100w60.004.txt & 764 & NA & 657 & - & 100.00 & 656 & 764 & 14.14 \\
                Dumas & n150w40.001.txt & 918 & NA & 810 & 919 & 11.86 & 812 & 918 & 11.55 \\
                Dumas & n150w40.002.txt & 941 & NA & 804 & 945 & 14.92 & 803 & 941 & 14.67 \\
                Dumas & n150w40.003.txt & 727 & NA & 528 & 737 & 28.36 & 525 & 727 & 27.79 \\
\hline
            \end{tabular}
    }
    \caption{\textbf{(Part 1 of 2)} TSPTW Results of PnB at 2048 on open problems: Seeded and Unseeded: LB = Lower Bound, UB = Upper Bound, OG = Optimality Gap.}
    \label{fig:tsptwopendata}
 \end{table}

  \begin{table}[!ht]
    \centering
    \resizebox{!}{6.65cm}{%
        \begin{tabular}{|ccc|c|ccc|ccc|}
            \multicolumn{3}{|c|}{Problem Information} & \multicolumn{1}{l|}{Baldacci et. al.   (2012)} & \multicolumn{3}{c|}{Single Thread: 2048} & \multicolumn{3}{c|}{Seeded Single Thread:   2048} \\
            Set & Name & Best Known Solution & LB & LB & UB & OG & LB & UB & OG \\ \hline
            AFG & rbg049a.tw & 10018 & 10006.10 & 9877 & 10034 & 1.56 & 9877 & 10018 & 1.41 \\
            AFG & rbg050b.tw & 9863 & 9855.80 & 9711 & 9876 & 1.67 & 9713 & 9863 & 1.52 \\
            AFG & rbg050c.tw & 10024 & 10020.40 & 9745 & 10074 & 3.27 & 9746 & 10024 & 2.77 \\
            AFG & rbg132.2.tw & 8192 & 8188.20 & 6701 & 8316 & 19.42 & 6695 & 8192 & 18.27 \\
            AFG & rbg152.3.tw & 9788 & 9785.20 & 7240 & 9982 & 27.47 & 7240 & 9788 & 26.03 \\
            AFG & rbg193.2.tw & 12142 & 12136.30 & 9755 & 12441 & 21.59 & 9755 & 12142 & 19.66 \\
\hline
        \end{tabular}
    }
    \caption{\textbf{(Part 2 of 2)} TSPTW Results of PnB at 2048 on open problems: Seeded and Unseeded: LB = Lower Bound, UB = Upper Bound, OG = Optimality Gap.}
    \label{fig:tsptwopendata2}
 \end{table}

\begin{table}[!ht]
    \centering
    \resizebox{!}{9.7cm}{%
        \begin{tabular}{|ccc|c|cccc|cccl}
            \multicolumn{3}{|c|}{Problem Information} & \multicolumn{1}{l|}{Baldacci et. al.   (2012)} & \multicolumn{4}{c|}{Single Thread: 2048} & \multicolumn{4}{c}{Seeded Single   Thread: 2048} \\
            Set & Name & Best Known Solution & LB & LB & UB & OG & Time & LB & UB & OG & \multicolumn{1}{c|}{Time} \\ \hline
            AFG & rbg010a.tw & 671 & O & 671 & 671 & 0 & 0.08 & 671 & 671 & 0 & \multicolumn{1}{l|}{0.26} \\
            AFG & rbg016a.tw & 938 & O & 938 & 938 & 0 & 0.22 & 938 & 938 & 0 & \multicolumn{1}{l|}{0.28} \\
            AFG & rbg016b.tw & 1304 & O & 1304 & 1304 & 0 & 5.62 & 1304 & 1304 & 0 & \multicolumn{1}{l|}{4.43} \\
            AFG & rbg017.2.tw & 852 & O & 852 & 852 & 0 & 9.08 & 852 & 852 & 0 & \multicolumn{1}{l|}{4.71} \\
            AFG & rbg017.tw & 893 & O & 893 & 893 & 0 & 3.92 & 893 & 893 & 0 & \multicolumn{1}{l|}{2.96} \\
            AFG & rbg017a.tw & 4296 & O & 4296 & 4296 & 0 & 0.28 & 4296 & 4296 & 0 & \multicolumn{1}{l|}{0.22} \\
\hline
            \end{tabular}
    }
    \caption{\textbf{(Part 1 of 6)} TSPTW Results of PnB at 2048 on closed problems: Seeded and Unseeded: LB = Lower Bound, UB = Upper Bound, OG = Optimality Gap, O = Optimality Proved.}
    \label{fig:tsptwcloseddata}
 \end{table}

 \begin{table}[!ht]
    \centering
    \resizebox{!}{9.7cm}{%
        \begin{tabular}{|ccc|c|cccc|cccl|}
            \multicolumn{3}{|c|}{Problem Information} & \multicolumn{1}{l|}{Baldacci et. al.   (2012)} & \multicolumn{4}{c|}{Single Thread: 2048} & \multicolumn{4}{c|}{Seeded   Single Thread: 2048} \\
            Set & Name & Best Known Solution & LB & LB & UB & OG & Time & LB & UB & OG & \multicolumn{1}{c|}{Time} \\ \hline
            Dumas & n20w20.001.txt & 378 & NA & 378 & 378 & 0 & 0.06 & 378 & 378 & 0 & 0.18 \\
            Dumas & n20w20.002.txt & 286 & NA & 286 & 286 & 0 & 0.08 & 286 & 286 & 0 & 0.18 \\
            Dumas & n20w20.003.txt & 394 & NA & 394 & 394 & 0 & 0.07 & 394 & 394 & 0 & 0.18 \\
            Dumas & n20w20.004.txt & 396 & NA & 396 & 396 & 0 & 0.06 & 396 & 396 & 0 & 0.18 \\
            Dumas & n20w20.005.txt & 352 & NA & 352 & 352 & 0 & 0.07 & 352 & 352 & 0 & 0.26 \\
            Dumas & n20w40.001.txt & 254 & NA & 254 & 254 & 0 & 0.12 & 254 & 254 & 0 & 0.31 \\
\hline
            \end{tabular}
    }
    \caption{\textbf{(Part 2 of 6)} TSPTW Results of PnB at 2048 on closed problems: Seeded and Unseeded: LB = Lower Bound, UB = Upper Bound, OG = Optimality Gap, O = Optimality Proved.}
    \label{fig:tsptwcloseddata2}
 \end{table}

 \begin{table}[!ht]
    \centering
    \resizebox{!}{9.7cm}{%
        \begin{tabular}{|ccc|c|cccc|cccl|}
            \multicolumn{3}{|c|}{Problem Information} & \multicolumn{1}{l|}{Baldacci et. al.   (2012)} & \multicolumn{4}{c|}{Single Thread: 2048} & \multicolumn{4}{c|}{Seeded   Single Thread: 2048} \\
            Set & Name & Best Known Solution & LB & LB & UB & OG & Time & LB & UB & OG & \multicolumn{1}{c|}{Time} \\ \hline
            Dumas & n60w80.005.txt & 468 & NA & 468 & 468 & 0 & 363.20 & 468 & 468 & 0 & 363.03 \\
            Dumas & n80w20.001.txt & 616 & NA & 616 & 616 & 0 & 58.66 & 616 & 616 & 0 & 31.87 \\
            Dumas & n80w20.002.txt & 737 & NA & 737 & 737 & 0 & 30.47 & 737 & 737 & 0 & 20.49 \\
            Dumas & n80w20.003.txt & 667 & NA & 667 & 667 & 0 & 14.48 & 667 & 667 & 0 & 5.52 \\
            Dumas & n80w20.004.txt & 615 & NA & 615 & 615 & 0 & 41.62 & 615 & 615 & 0 & 32.25 \\
            Dumas & n80w20.005.txt & 748 & NA & 748 & 748 & 0 & 23.05 & 748 & 748 & 0 & 10.52 \\
\hline
            \end{tabular}
    }
    \caption{\textbf{(Part 3 of 6)} TSPTW Results of PnB at 2048 on closed problems: Seeded and Unseeded: LB = Lower Bound, UB = Upper Bound, OG = Optimality Gap, O = Optimality Proved.}
    \label{fig:tsptwcloseddata3}
 \end{table}

 \begin{table}[!ht]
    \centering
    \resizebox{!}{9.7cm}{%
        \begin{tabular}{|ccc|c|cccc|cccl|}
            \multicolumn{3}{|c|}{Problem Information} & \multicolumn{1}{l|}{Baldacci et. al.   (2012)} & \multicolumn{4}{c|}{Single Thread: 2048} & \multicolumn{4}{c|}{Seeded   Single Thread: 2048} \\
            Set & Name & Best Known Solution & LB & LB & UB & OG & Time & LB & UB & OG & \multicolumn{1}{c|}{Time} \\ \hline
            GendreauDumas & n40w140.004.txt & 342 & O & 246 & 342 & 28.07 & - & 246 & 342 & 28.07 & - \\
            GendreauDumas & n40w140.005.txt & 371 & O & 371 & 371 & 0 & 471.26 & 371 & 371 & 0 & 289.28 \\
            GendreauDumas & n40w160.001.txt & 348 & O & 348 & 348 & 0 & 179.81 & 348 & 348 & 0 & 159.90 \\
            GendreauDumas & n40w160.002.txt & 337 & O & 337 & 337 & 0 & 127.84 & 337 & 337 & 0 & 115.75 \\
            GendreauDumas & n40w160.003.txt & 346 & O & 346 & 346 & 0 & 301.12 & 346 & 346 & 0 & 156.85 \\
            GendreauDumas & n40w160.004.txt & 288 & O & 288 & 288 & 0 & 313.23 & 288 & 288 & 0 & 292.72 \\
\hline
            \end{tabular}
    }
    \caption{\textbf{(Part 4 of 6)} TSPTW Results of PnB at 2048 on closed problems: Seeded and Unseeded: LB = Lower Bound, UB = Upper Bound, OG = Optimality Gap, O = Optimality Proved.}
    \label{fig:tsptwcloseddata4}
 \end{table}

 \begin{table}[!ht]
    \centering
    \resizebox{!}{9.7cm}{%
        \begin{tabular}{|ccc|c|cccc|cccl|}
            \multicolumn{3}{|c|}{Problem Information} & \multicolumn{1}{l|}{Baldacci et. al.   (2012)} & \multicolumn{4}{c|}{Single Thread: 2048} & \multicolumn{4}{c|}{Seeded   Single Thread: 2048} \\
            Set & Name & Best Known Solution & LB & LB & UB & OG & Time & LB & UB & OG & \multicolumn{1}{c|}{Time} \\ \hline
            Langevin & N40ft208.dat & 1050.7 & O & 1050.7 & 1050.7 & 0 & 0.08 & 1050.7 & 1050.7 & 0 & 0.24 \\
            Langevin & N40ft209.dat & 1013.9 & O & 1013.9 & 1013.9 & 0 & 0.07 & 1013.9 & 1013.9 & 0 & 0.29 \\
            Langevin & N40ft210.dat & 1026.3 & O & 1026.3 & 1026.3 & 0 & 0.08 & 1026.3 & 1026.3 & 0 & 0.28 \\
            Langevin & N40ft401.dat & 1085 & O & 1085 & 1085 & 0 & 0.11 & 1085 & 1085 & 0 & 0.79 \\
            Langevin & N40ft402.dat & 995.6 & O & 995.6 & 995.6 & 0 & 0.07 & 995.6 & 995.6 & 0 & 0.34 \\
            Langevin & N40ft403.dat & 845.8 & O & 845.8 & 845.8 & 0 & 2.03 & 845.8 & 845.8 & 0 & 0.41 \\
\hline
            \end{tabular}
    }
    \caption{\textbf{(Part 5 of 6)} TSPTW Results of PnB at 2048 on closed problems: Seeded and Unseeded: LB = Lower Bound, UB = Upper Bound, OG = Optimality Gap, O = Optimality Proved.}
    \label{fig:tsptwcloseddata5}
 \end{table}

 \begin{table}[!ht]
    \centering
    \resizebox{!}{5.82cm}{%
        \begin{tabular}{|ccc|c|cccc|cccl|}
            \multicolumn{3}{|c|}{Problem Information} & \multicolumn{1}{l|}{Baldacci et. al.   (2012)} & \multicolumn{4}{c|}{Single Thread: 2048} & \multicolumn{4}{c|}{Seeded   Single Thread: 2048} \\
            Set & Name & Best Known Solution & LB & LB & UB & OG & Time & LB & UB & OG & \multicolumn{1}{c|}{Time} \\ \hline
            SolomonPesant & rc201.0 & 628.62 & O & 628.62 & 628.62 & 0 & 3.36 & 628.62 & 628.62 & 0 & 2.13 \\
            SolomonPesant & rc201.1 & 654.7 & O & 654.70 & 654.70 & 0 & 3.97 & 654.70 & 654.70 & 0 & 2.56 \\
            SolomonPesant & rc201.2 & 707.65 & O & 707.65 & 707.65 & 0 & 3.34 & 707.65 & 707.65 & 0 & 1.63 \\
            SolomonPesant & rc201.3 & 422.54 & O & 422.54 & 422.54 & 0 & 0.52 & 422.54 & 422.54 & 0 & 0.73 \\
            SolomonPesant & rc202.0 & 496.22 & O & 496.22 & 496.22 & 0 & 19.91 & 496.22 & 496.22 & 0 & 12.48 \\
            SolomonPesant & rc202.1 & 426.53 & O & 426.53 & 426.53 & 0 & 14.76 & 426.53 & 426.53 & 0 & 10.34 \\
\hline
            \end{tabular}
    }
    \caption{\textbf{(Part 6 of 6)} TSPTW Results of PnB at 2048 on closed problems: Seeded and Unseeded: LB = Lower Bound, UB = Upper Bound, OG = Optimality Gap, O = Optimality Proved.}
    \label{fig:tsptwcloseddata6}
 \end{table}

 \begin{table}[!ht]
    \centering
    \resizebox{!}{9.7cm}{%
        \begin{tabular}{|ccc|ccc|ccc|}
            \multicolumn{3}{|c|}{Problem Information} & \multicolumn{3}{c|}{Single Thread: 2048} & \multicolumn{3}{c|}{Seeded   Single Thread: 2048} \\
            Set & Name & Best Known Solution & LB & UB & Time & LB & UB & Time \\ \hline
            AFG & rbg010a.tw & 3840 & 3840 & 3840 & 0.10 & 3840 & 3840 & 0.27 \\
            AFG & rbg016a.tw & 2596 & 2596 & 2596 & 0.59 & 2596 & 2596 & 0.23 \\
            AFG & rbg016b.tw & 2094 & 2094 & 2094 & 1.79 & 2094 & 2094 & 0.74 \\
            AFG & rbg017.2.tw & 2351 & 2351 & 2351 & 4.68 & 2351 & 2351 & 2.93 \\
            AFG & rbg017.tw & 2351 & 2351 & 2351 & 1.55 & 2351 & 2351 & 0.97 \\
            AFG & rbg017a.tw & 4296 & 4296 & 4296 & 0.29 & 4296 & 4296 & 0.25 \\
\hline
            \end{tabular}
    }
    \caption{\textbf{(Part 1 of 7)} Makespan Results of PnB at 2048 on closed problems: Seeded and Unseeded: LB = Lower Bound, UB = Upper Bound.}
    \label{fig:makespanfulldata}
 \end{table}

 \begin{table}[!ht]
    \centering
    \resizebox{!}{9.7cm}{%
        \begin{tabular}{|ccc|ccc|ccc|}
            \multicolumn{3}{|c|}{Problem Information} & \multicolumn{3}{c|}{Single Thread: 2048} & \multicolumn{3}{c|}{Seeded   Single Thread: 2048} \\
            Set & Name & Best Known Solution & LB & UB & Time & LB & UB & Time \\ \hline
            Dumas & n150w40.001.txt & 1058 & 1058 & 1058 & 1189.38 & 1058 & 1058 & 183.42 \\
            Dumas & n150w40.002.txt & 1072 & 1072 & 1072 & 302.80 & 1072 & 1072 & 227.35 \\
            Dumas & n150w40.003.txt & 894 & 894 & 894 & 202.91 & 894 & 894 & 495.32 \\
            Dumas & n150w40.004.txt & 948 & 948 & 948 & 261.44 & 948 & 948 & 162.54 \\
            Dumas & n150w40.005.txt & 980 & 980 & 980 & 3533.27 & 980 & 980 & 194.03 \\
            Dumas & n150w60.001.txt & 1024 & 1024 & 1024 & 1306.06 & 1024 & 1024 & 1361.14 \\
\hline
            \end{tabular}
    }
    \caption{\textbf{(Part 2 of 7)} Makespan Results of PnB at 2048 on closed problems: Seeded and Unseeded: LB = Lower Bound, UB = Upper Bound.}
    \label{fig:makespanfulldata2}
 \end{table}

  \begin{table}[!ht]
    \centering
    \resizebox{!}{9.7cm}{%
        \begin{tabular}{|ccc|ccc|ccc|}
            \multicolumn{3}{|c|}{Problem Information} & \multicolumn{3}{c|}{Single Thread: 2048} & \multicolumn{3}{c|}{Seeded   Single Thread: 2048} \\
            Set & Name & Best Known Solution & LB & UB & Time & LB & UB & Time \\ \hline
            Dumas & n60w100.001.txt & 576 & 576 & 576 & 47.90 & 576 & 576 & 41.17 \\
            Dumas & n60w100.002.txt & 622 & 622 & 622 & 49.11 & 622 & 622 & 31.61 \\
            Dumas & n60w100.003.txt & 610 & 610 & 610 & 51.91 & 610 & 610 & 28.48 \\
            Dumas & n60w100.004.txt & 625 & 625 & 625 & 56.84 & 625 & 625 & 55.01 \\
            Dumas & n60w100.005.txt & 668 & 668 & 668 & 85.20 & 668 & 668 & 42.20 \\
            Dumas & n60w20.001.txt & 586 & 586 & 586 & 6.05 & 586 & 586 & 2.00 \\
\hline
            \end{tabular}
    }
    \caption{\textbf{(Part 3 of 7)} Makespan Results of PnB at 2048 on closed problems: Seeded and Unseeded: LB = Lower Bound, UB = Upper Bound.}
    \label{fig:makespanfulldata3}
 \end{table}
 
 \begin{table}[!ht]
    \centering
    \resizebox{!}{9.7cm}{%
        \begin{tabular}{|ccc|ccc|ccc|}
            \multicolumn{3}{|c|}{Problem Information} & \multicolumn{3}{c|}{Single Thread: 2048} & \multicolumn{3}{c|}{Seeded   Single Thread: 2048} \\
            Set & Name & Best Known Solution & LB & UB & Time & LB & UB & Time \\ \hline
            GendreauDumas & n20w120.001.txt & 337 & 337 & 337 & 5.13 & 337 & 337 & 4.10 \\
            GendreauDumas & n20w120.002.txt & 246 & 246 & 246 & 6.13 & 246 & 246 & 5.79 \\
            GendreauDumas & n20w120.003.txt & 347 & 347 & 347 & 8.12 & 347 & 347 & 3.19 \\
            GendreauDumas & n20w120.004.txt & 353 & 353 & 353 & 9.36 & 353 & 353 & 5.69 \\
            GendreauDumas & n20w120.005.txt & 315 & 315 & 315 & 11.88 & 315 & 315 & 7.86 \\
            GendreauDumas & n20w140.001.txt & 230 & 230 & 230 & 7.33 & 230 & 230 & 5.63 \\
\hline
        \end{tabular}
    }
    \caption{\textbf{(Part 4 of 7)} Makespan Results of PnB at 2048 on closed problems: Seeded and Unseeded: LB = Lower Bound, UB = Upper Bound.}
    \label{fig:makespanfulldata4}
 \end{table}

  \begin{table}[!ht]
    \centering
    \resizebox{!}{9.7cm}{%
        \begin{tabular}{|ccc|ccc|ccc|}
            \multicolumn{3}{|c|}{Problem Information} & \multicolumn{3}{c|}{Single Thread: 2048} & \multicolumn{3}{c|}{Seeded   Single Thread: 2048} \\
            Set & Name & Best Known Solution & LB & UB & Time & LB & UB & Time \\ \hline
            GendreauDumas & n60w180.004.txt & 629 & 629 & 629 & 264.17 & 629 & 629 & 110.83 \\
            GendreauDumas & n60w180.005.txt & 528 & 528 & 528 & 219.24 & 528 & 528 & 102.20 \\
            GendreauDumas & n60w200.001.txt & 526 & 526 & 526 & 846.26 & 526 & 526 & 645.07 \\
            GendreauDumas & n60w200.002.txt & 572 & 572 & 572 & 235.71 & 572 & 572 & 68.66 \\
            GendreauDumas & n60w200.004.txt & 575 & 575 & 575 & 286.27 & 575 & 575 & 53.71 \\
            GendreauDumas & n60w200.005.txt & 618 & 618 & 618 & 292.92 & 618 & 618 & 159.49 \\
\hline
            \end{tabular}
    }
    \caption{\textbf{(Part 5 of 7)} Makespan Results of PnB at 2048 on closed problems: Seeded and Unseeded: LB = Lower Bound, UB = Upper Bound.}
    \label{fig:makespanfulldata5}
 \end{table}

 \begin{table}[!ht]
    \centering
    \resizebox{!}{9.7cm}{%
        \begin{tabular}{|ccc|ccc|ccc|}
            \multicolumn{3}{|c|}{Problem Information} & \multicolumn{3}{c|}{Single Thread: 2048} & \multicolumn{3}{c|}{Seeded   Single Thread: 2048} \\
            Set & Name & Best Known Solution & LB & UB & Time & LB & UB & Time \\ \hline
            Langevin & N40ft408.dat & 1043.5 & 1043.5 & 1043.5 & 0.25 & 1043.5 & 1043.5 & 0.33 \\
            Langevin & N40ft409.dat & 1025.5 & 1025.5 & 1025.5 & 0.12 & 1025.5 & 1025.5 & 0.35 \\
            Langevin & N40ft410.dat & 1034.3 & 1034.3 & 1034.3 & 0.07 & 1034.3 & 1034.3 & 0.30 \\
            Langevin & N60ft201.dat & 1375.4 & 1375.4 & 1375.4 & 0.13 & 1375.4 & 1375.4 & 0.46 \\
            Langevin & N60ft202.dat & 1186.4 & 1186.4 & 1186.4 & 0.09 & 1186.4 & 1186.4 & 0.44 \\
            Langevin & N60ft203.dat & 1194.2 & 1194.2 & 1194.2 & 0.08 & 1194.2 & 1194.2 & 0.45 \\
\hline
            \end{tabular}
    }
    \caption{\textbf{(Part 6 of 7)} Makespan Results of PnB at 2048 on closed problems: Seeded and Unseeded: LB = Lower Bound, UB = Upper Bound.}
    \label{fig:makespanfulldata6}
 \end{table}

 \begin{table}[!b]
    \centering
    \resizebox{!}{4.85cm}{%
        \begin{tabular}{|ccc|ccc|ccc|}
            \multicolumn{3}{|c|}{Problem Information} & \multicolumn{3}{c|}{Single Thread: 2048} & \multicolumn{3}{c|}{Seeded   Single Thread: 2048} \\
            Set & Name & Best Known Solution & LB & UB & Time & LB & UB & Time \\ \hline
            SolomonPesant & rc207.1 & 785.37 & 785.37 & 785.37 & 132.40 & 785.37 & 785.37 & 39.23 \\
            SolomonPesant & rc207.2 & 650.8 & 650.80 & 650.80 & 252.09 & 650.80 & 650.80 & 35.39 \\
            SolomonPesant & rc208.0 & 836.04 & 836.04 & 836.04 & 2256.13 & 836.04 & 836.04 & 1484.62 \\
            SolomonPesant & rc208.1 & 615.51 & 615.51 & 615.51 & 417.09 & 615.51 & 615.51 & 160.41 \\
            SolomonPesant & rc208.2 & 596.21 & 596.21 & 596.21 & 276.60 & 596.21 & 596.21 & 46.52 \\ \hline
            SolomonPotvinBengio & rc\_201.1.txt & 592.06 & 592.06 & 592.06 & 1.75 & 592.06 & 592.06 & 0.53 \\
\hline
            \end{tabular}
    }
    \caption{\textbf{(Part 7 of 7)} Makespan Results of PnB at 2048 on closed problems: Seeded and Unseeded: LB = Lower Bound, UB = Upper Bound.}
    \label{fig:makespanfulldata7}
 \end{table}
\FloatBarrier
\bibliography{peel-and-bound}

\begin{thebibliography}{}

\bibitem[\protect\BCAY{Andersen, Hadzic, Hooker,\ \BBA\ Tiedemann}{Andersen
  et~al.}{2007}]{andersen2007}
Andersen, H., Hadzic, T., Hooker, J., \BBA\ Tiedemann, P. \BBOP2007\BBCP.
\newblock \BBOQ A constraint store based on multivalued decision diagrams\BBCQ\
\newblock In {\Bem Bessière, C. (eds) Principles and Practice of Constraint
  Programming – CP 2007}, \lowercase{\BVOL}\ 4741 of {\Bem Lecture Notes in
  Computer Science}, \BPGS\ 118--132.

\bibitem[\protect\BCAY{Ascheuer}{Ascheuer}{1996}]{AFG}
Ascheuer, N. \BBOP1996\BBCP.
\newblock {\Bem Hamiltonian path problems in the on-line optimization of
  flexible manufacturing systems}.
\newblock Ph.D.\ thesis, Konrad–Zuse–Zentrum f\"ur Informationstechnik
  Berlin.

\bibitem[\protect\BCAY{Baldacci, Mingozzi,\ \BBA\ Roberti}{Baldacci
  et~al.}{2012}]{baldacci2012new}
Baldacci, R., Mingozzi, A., \BBA\ Roberti, R. \BBOP2012\BBCP.
\newblock \BBOQ New state-space relaxations for solving the traveling salesman
  problem with time windows\BBCQ\
\newblock {\Bem INFORMS Journal on Computing}, {\Bem 24\/}(3), 356--371.

\bibitem[\protect\BCAY{Baptiste, Le~Pape,\ \BBA\ Nuijten}{Baptiste
  et~al.}{2001}]{cpS}
Baptiste, P., Le~Pape, C., \BBA\ Nuijten, W. \BBOP2001\BBCP.
\newblock {\Bem Constraint-Based Scheduling: Applying Constraint Programming to
  Scheduling Problems}.
\newblock International Series in Operations Research and Management Science,
  Kluwer. Springer International Publishing.

\bibitem[\protect\BCAY{Bergman, Cire, van Hoeve,\ \BBA\ Hooker}{Bergman
  et~al.}{2016a}]{DDforO}
Bergman, D., Cire, A., van Hoeve, W.-J., \BBA\ Hooker, J. \BBOP2016a\BBCP.
\newblock {\Bem Decision Diagrams for Optimization}.
\newblock Springer International Publishing.

\bibitem[\protect\BCAY{Bergman, Cire, van Hoeve,\ \BBA\ Hooker}{Bergman
  et~al.}{2016b}]{BnB2}
Bergman, D., Cire, A., van Hoeve, W.-J., \BBA\ Hooker, J. \BBOP2016b\BBCP.
\newblock \BBOQ Discrete optimization with decision diagrams\BBCQ\
\newblock {\Bem INFORMS Journal on Computing}, {\Bem 28}, 47--66.

\bibitem[\protect\BCAY{Bergman, Cire, van Hoeve,\ \BBA\ Hooker}{Bergman
  et~al.}{2014a}]{bergman2014optimization}
Bergman, D., Cire, A.~A., van Hoeve, W.-J., \BBA\ Hooker, J.~N.
  \BBOP2014a\BBCP.
\newblock \BBOQ Optimization bounds from binary decision diagrams\BBCQ\
\newblock {\Bem INFORMS Journal on Computing}, {\Bem 26\/}(2), 253--268.

\bibitem[\protect\BCAY{Bergman, Cire, Sabharwal, Samulowitz, Saraswat,\ \BBA\
  van Hoeve}{Bergman et~al.}{2014b}]{BnB}
Bergman, D., Cire, A.~A., Sabharwal, A., Samulowitz, H., Saraswat, V., \BBA\
  van Hoeve, W.-J. \BBOP2014b\BBCP.
\newblock \BBOQ Parallel combinatorial optimization with decision
  diagrams\BBCQ\
\newblock In {\Bem Proceedings of the International Conference on AI and OR
  Techniques in Constraint Programming for Combinatorial Optimization
  Problems}, \BPGS\ 351--367.

\bibitem[\protect\BCAY{Bergman, van Hoeve,\ \BBA\ Hooker}{Bergman
  et~al.}{2011}]{bergman2011}
Bergman, D., van Hoeve, W.-J., \BBA\ Hooker, J. \BBOP2011\BBCP.
\newblock \BBOQ Manipulating mdd relaxations for combinatorial
  optimization\BBCQ\
\newblock In {\Bem Lecture Notes in Computer Science}, \lowercase{\BVOL}\ 6697,
  \BPGS\ 20--35.

\bibitem[\protect\BCAY{Cappart, Bergman, Rousseau, Pr{\'e}mont-Schwarz,\ \BBA\
  Parjadis}{Cappart et~al.}{2022}]{cappart2022improving}
Cappart, Q., Bergman, D., Rousseau, L.-M., Pr{\'e}mont-Schwarz, I., \BBA\
  Parjadis, A. \BBOP2022\BBCP.
\newblock \BBOQ Improving variable orderings of approximate decision diagrams
  using reinforcement learning\BBCQ\
\newblock {\Bem INFORMS Journal on Computing}, {\Bem 34\/}(5).

\bibitem[\protect\BCAY{Cappart, Goutierre, Bergman,\ \BBA\ Rousseau}{Cappart
  et~al.}{2019}]{cappart2019improving}
Cappart, Q., Goutierre, E., Bergman, D., \BBA\ Rousseau, L.-M. \BBOP2019\BBCP.
\newblock \BBOQ Improving optimization bounds using machine learning: Decision
  diagrams meet deep reinforcement learning\BBCQ\
\newblock In {\Bem Proceedings of the AAAI Conference on Artificial
  Intelligence}, \lowercase{\BVOL}~33, \BPGS\ 1443--1451.

\bibitem[\protect\BCAY{Castro, Piacentini, Cire,\ \BBA\ Beck}{Castro
  et~al.}{2020}]{dual2}
Castro, M., Piacentini, C., Cire, A., \BBA\ Beck, J. \BBOP2020\BBCP.
\newblock \BBOQ Solving delete free planning with relaxed decision diagram
  based heuristics\BBCQ\
\newblock {\Bem Journal of Artificial Intelligence Research}, {\Bem 67},
  607--651.

\bibitem[\protect\BCAY{Castro, Cire,\ \BBA\ Beck}{Castro
  et~al.}{2022}]{ddsurvey}
Castro, M.~P., Cire, A.~A., \BBA\ Beck, J.~C. \BBOP2022\BBCP.
\newblock \BBOQ Decision diagrams for discrete optimization: A survey of recent
  advances\BBCQ\
\newblock {\Bem INFORMS Journal on Computing}, {\Bem 34}.

\bibitem[\protect\BCAY{Cire\ \BBA\ van Hoeve}{Cire\ \BBA\ van
  Hoeve}{2013}]{MDDForSP}
Cire, A.~A.\BBACOMMA\  \BBA\ van Hoeve, W.-J. \BBOP2013\BBCP.
\newblock \BBOQ Multivalued decision diagrams for sequencing problems\BBCQ\
\newblock {\Bem Operations Research}, {\Bem 61\/}(6), 1259, 1462.

\bibitem[\protect\BCAY{Copp'e, Gillard,\ \BBA\ Schaus}{Copp'e
  et~al.}{2022}]{vianney}
Copp'e, V., Gillard, X., \BBA\ Schaus, P. \BBOP2022\BBCP.
\newblock \BBOQ Branch-and-bound with barrier: Dominance and suboptimality
  detection for dd-based branch-and-bound\BBCQ\
\newblock {\Bem ArXiv}, {\Bem abs/2211.13118}.

\bibitem[\protect\BCAY{Dumas, Desrosiers, Gelinas,\ \BBA\ Solomon}{Dumas
  et~al.}{1995}]{Dumas}
Dumas, Y., Desrosiers, J., Gelinas, E., \BBA\ Solomon, M.~M. \BBOP1995\BBCP.
\newblock \BBOQ An optimal algorithm for the traveling salesman problem with
  time windows\BBCQ\
\newblock {\Bem Operations research}, {\Bem 43\/}(2), 367--371.

\bibitem[\protect\BCAY{Gendreau, Hertz, Laporte,\ \BBA\ Stan}{Gendreau
  et~al.}{1998}]{gendreaudumas}
Gendreau, M., Hertz, A., Laporte, G., \BBA\ Stan, M. \BBOP1998\BBCP.
\newblock \BBOQ A generalized insertion heuristic for the traveling salesman
  problem with time windows\BBCQ\
\newblock {\Bem Operations Research}, {\Bem 46\/}(3), 330--335.

\bibitem[\protect\BCAY{Gentzel, Michel,\ \BBA\ Hoeve}{Gentzel
  et~al.}{2020}]{gentzel2020haddock}
Gentzel, R., Michel, L., \BBA\ Hoeve, W.-J.~v. \BBOP2020\BBCP.
\newblock \BBOQ Haddock: A language and architecture for decision diagram
  compilation\BBCQ\
\newblock In {\Bem International Conference on Principles and Practice of
  Constraint Programming}, \BPGS\ 531--547. Springer.

\bibitem[\protect\BCAY{Gillard}{Gillard}{2022}]{gillard2022discrete}
Gillard, X. \BBOP2022\BBCP.
\newblock {\Bem Discrete optimization with decision diagrams: design of a
  generic solver, improved bounding techniques, and discovery of good feasible
  solutions with large neighborhood search}.
\newblock Ph.D.\ thesis, UCL-Universit{\'e} Catholique de Louvain.

\bibitem[\protect\BCAY{Gillard, Coppé, Schaus,\ \BBA\ Cire}{Gillard
  et~al.}{2021}]{ddo}
Gillard, X., Coppé, V., Schaus, P., \BBA\ Cire, A.~A. \BBOP2021\BBCP.
\newblock \BBOQ Improving the filtering of branch-and-bound mdd solver\BBCQ\
\newblock In {\Bem Integration of Constraint Programming, Artificial
  Intelligence, and Operations Research, 18th International Conference, CPAIOR
  2021}, Lecture Notes in Computer Science. Springer International Publishing.

\bibitem[\protect\BCAY{González, Cire, Lodi,\ \BBA\ Rousseau}{González
  et~al.}{2020}]{bnb3}
González, J., Cire, A., Lodi, A., \BBA\ Rousseau, L.-M. \BBOP2020\BBCP.
\newblock \BBOQ Integrated integer programming and decision diagram search tree
  with an application to the maximum independent set problem\BBCQ\
\newblock {\Bem Constraints}, {\Bem 25}.

\bibitem[\protect\BCAY{Hadzic, Hooker, O'Sullivan,\ \BBA\ Tiedemann}{Hadzic
  et~al.}{2008}]{hadzic2008}
Hadzic, T., Hooker, J., O'Sullivan, B., \BBA\ Tiedemann, P. \BBOP2008\BBCP.
\newblock \BBOQ Approximate compilation of constraints into multivalued
  decision diagrams\BBCQ\
\newblock In {\Bem Stuckey, P.J. (eds) Principles and Practice of Constraint
  Programming. CP 2008}, Lecture Notes in Computer Science, \BPGS\ 448--462.

\bibitem[\protect\BCAY{Hoda, van Hoeve,\ \BBA\ Hooker}{Hoda
  et~al.}{2010}]{hoda2010}
Hoda, S., van Hoeve, W.-J., \BBA\ Hooker, J. \BBOP2010\BBCP.
\newblock \BBOQ A systematic approach to mdd-based constraint programming\BBCQ\
\newblock In {\Bem Cohen, D. (eds) Principles and Practice of Constraint
  Programming – CP 2010. CP 2010}, \lowercase{\BVOL}\ 6308 of {\Bem Lecture
  Notes in Computer Science}, \BPGS\ 266--280.

\bibitem[\protect\BCAY{Hoeve}{Hoeve}{2021}]{dual1}
Hoeve, W.-J. \BBOP2021\BBCP.
\newblock \BBOQ Graph coloring with decision diagrams\BBCQ\
\newblock {\Bem Mathematical Programming}, {\Bem 192}, 631--674.

\bibitem[\protect\BCAY{Hooker}{Hooker}{2013}]{canonicalarcs}
Hooker, J. \BBOP2013\BBCP.
\newblock \BBOQ Decision diagrams and dynamic programming\BBCQ\
\newblock In {\Bem Gomes, C., Sellmann, M. (eds) Integration of AI and OR
  Techniques in Constraint Programming for Combinatorial Optimization Problems.
  CPAIOR 2013}, \lowercase{\BVOL}\ 7874 of {\Bem Lecture Notes in Computer
  Science}.

\bibitem[\protect\BCAY{Hooker}{Hooker}{2017}]{jobsequencing}
Hooker, J. \BBOP2017\BBCP.
\newblock \BBOQ Job sequencing bounds from decision diagrams\BBCQ\
\newblock In {\Bem Beck, J. (eds) Principles and Practice of Constraint
  Programming. CP 2017}, Lecture Notes in Computer Science, \BPGS\ 565--578.

\bibitem[\protect\BCAY{Hooker}{Hooker}{2019}]{jobsequencing2}
Hooker, J. \BBOP2019\BBCP.
\newblock \BBOQ Improved job sequencing bounds from decision diagrams\BBCQ\
\newblock In {\Bem Schiex, T., de Givry, S. (eds) Principles and Practice of
  Constraint Programming. CP 2019}, \lowercase{\BVOL}\ 11802 of {\Bem Lecture
  Notes in Computer Science}, \BPGS\ 268--283.

\bibitem[\protect\BCAY{Karahalios\ \BBA\ Hoeve}{Karahalios\ \BBA\
  Hoeve}{2022}]{karahalios}
Karahalios, A.\BBACOMMA\  \BBA\ Hoeve, W.-J. \BBOP2022\BBCP.
\newblock \BBOQ Variable ordering for decision diagrams: A portfolio
  approach\BBCQ\
\newblock {\Bem Constraints}, {\Bem 27\/}(1-2), 1--18.

\bibitem[\protect\BCAY{Langevin, Desrochers, Desrosiers, G{\'e}linas,\ \BBA\
  Soumis}{Langevin et~al.}{1993}]{langevin}
Langevin, A., Desrochers, M., Desrosiers, J., G{\'e}linas, S., \BBA\ Soumis, F.
  \BBOP1993\BBCP.
\newblock \BBOQ A two-commodity flow formulation for the traveling salesman and
  the makespan problems with time windows\BBCQ\
\newblock {\Bem Networks}, {\Bem 23\/}(7), 631--640.

\bibitem[\protect\BCAY{Latour, Babaki,\ \BBA\ Nijssen}{Latour
  et~al.}{2019}]{stochastic2}
Latour, A., Babaki, B., \BBA\ Nijssen, S. \BBOP2019\BBCP.
\newblock \BBOQ Stochastic constraint propagation for mining probabilistic
  networks\BBCQ\
\newblock In {\Bem Proceedings of the Twenty-Eighth International Joint
  Conference on Artificial Intelligence}, \BPGS\ 1137--1145.

\bibitem[\protect\BCAY{Lozano\ \BBA\ Smith}{Lozano\ \BBA\
  Smith}{2018}]{stochastic1}
Lozano, L.\BBACOMMA\  \BBA\ Smith, J. \BBOP2018\BBCP.
\newblock \BBOQ A binary decision diagram based algorithm for solving a class
  of binary two-stage stochastic programs\BBCQ\
\newblock {\Bem Mathematical Programming}, {\Bem 191}, 381--404.

\bibitem[\protect\BCAY{López-Ibáñez\ \BBA\ Blum}{López-Ibáñez\ \BBA\
  Blum}{2022}]{tsptwbench}
López-Ibáñez, M.\BBACOMMA\  \BBA\ Blum, C. \BBOP2022\BBCP.
\newblock \BBOQ Benchmark instances for the travelling salesman problem with
  time windows (tsptw)\BBCQ\
\newblock https://lopez-ibanez.eu/tsptw-instances.
\newblock Accessed: 2022-01-22.

\bibitem[\protect\BCAY{Maschler\ \BBA\ Raidl}{Maschler\ \BBA\
  Raidl}{2021}]{dual3}
Maschler, J.\BBACOMMA\  \BBA\ Raidl, G. \BBOP2021\BBCP.
\newblock \BBOQ Multivalued decision diagrams for prize-collecting job
  sequencing with one common and multiple secondary resources\BBCQ\
\newblock {\Bem Annals of Operations Research}, {\Bem 302}.

\bibitem[\protect\BCAY{Ohlmann\ \BBA\ Thomas}{Ohlmann\ \BBA\
  Thomas}{2007}]{ohlmann}
Ohlmann, J.~W.\BBACOMMA\  \BBA\ Thomas, B.~W. \BBOP2007\BBCP.
\newblock \BBOQ A compressed-annealing heuristic for the traveling salesman
  problem with time windows\BBCQ\
\newblock {\Bem INFORMS Journal on Computing}, {\Bem 19\/}(1), 80--90.

\bibitem[\protect\BCAY{Parjadis, Cappart, Rousseau,\ \BBA\ Bergman}{Parjadis
  et~al.}{2021}]{parjadis2021improving}
Parjadis, A., Cappart, Q., Rousseau, L.-M., \BBA\ Bergman, D. \BBOP2021\BBCP.
\newblock \BBOQ Improving branch-and-bound using decision diagrams and
  reinforcement learning\BBCQ\
\newblock In {\Bem International Conference on Integration of Constraint
  Programming, Artificial Intelligence, and Operations Research}, \BPGS\
  446--455. Springer.

\bibitem[\protect\BCAY{Perez\ \BBA\ Régin}{Perez\ \BBA\
  Régin}{2018}]{Perez_Regin_2018}
Perez, G.\BBACOMMA\  \BBA\ Régin, J.-C. \BBOP2018\BBCP.
\newblock \BBOQ Parallel algorithms for operations on multi-valued decision
  diagrams\BBCQ\
\newblock {\Bem Proceedings of the AAAI Conference on Artificial Intelligence},
  {\Bem 32\/}(1).

\bibitem[\protect\BCAY{Pesant, Gendreau, Potvin,\ \BBA\ Rousseau}{Pesant
  et~al.}{1998}]{pesant1998exact}
Pesant, G., Gendreau, M., Potvin, J.-Y., \BBA\ Rousseau, J.-M. \BBOP1998\BBCP.
\newblock \BBOQ An exact constraint logic programming algorithm for the
  traveling salesman problem with time windows\BBCQ\
\newblock {\Bem Transportation Science}, {\Bem 32\/}(1), 12--29.

\bibitem[\protect\BCAY{Potvin\ \BBA\ Bengio}{Potvin\ \BBA\
  Bengio}{1996}]{potvin1996vehicle}
Potvin, J.-Y.\BBACOMMA\  \BBA\ Bengio, S. \BBOP1996\BBCP.
\newblock \BBOQ The vehicle routing problem with time windows part ii: genetic
  search\BBCQ\
\newblock {\Bem INFORMS journal on Computing}, {\Bem 8\/}(2), 165--172.

\bibitem[\protect\BCAY{Reinelt}{Reinelt}{1991}]{tsplib}
Reinelt, G. \BBOP1991\BBCP.
\newblock \BBOQ Tsplib. a traveling salesman problem library\BBCQ\
\newblock {\Bem INFORMS Journal on Computing}, {\Bem 3}, 376--384.

\bibitem[\protect\BCAY{Rudich, Cappart,\ \BBA\ Rousseau}{Rudich
  et~al.}{2022}]{Rudich2022PeelAndBoundGS}
Rudich, I., Cappart, Q., \BBA\ Rousseau, L.-M. \BBOP2022\BBCP.
\newblock \BBOQ Peel-and-bound: Generating stronger relaxed bounds with
  multivalued decision diagrams\BBCQ\
\newblock In {\Bem International Conference on Principles and Practice of
  Constraint Programming}.

\bibitem[\protect\BCAY{Serra\ \BBA\ Hooker}{Serra\ \BBA\
  Hooker}{2019}]{postopt}
Serra, T.\BBACOMMA\  \BBA\ Hooker, J. \BBOP2019\BBCP.
\newblock \BBOQ Compact representation of near-optimal integer programming
  solutions\BBCQ\
\newblock {\Bem Mathematical Programming}, {\Bem 182}.

\bibitem[\protect\BCAY{Uña, Gange, Schachte,\ \BBA\ Stuckey}{Uña
  et~al.}{2019}]{global1}
Uña, D., Gange, G., Schachte, P., \BBA\ Stuckey, P. \BBOP2019\BBCP.
\newblock \BBOQ Compiling cp subproblems to mdds and d-dnnfs\BBCQ\
\newblock {\Bem Constraints}, {\Bem 24}.

\bibitem[\protect\BCAY{Verhaeghe, Lecoutre,\ \BBA\ Schaus}{Verhaeghe
  et~al.}{2018}]{global2}
Verhaeghe, H., Lecoutre, C., \BBA\ Schaus, P. \BBOP2018\BBCP.
\newblock \BBOQ Compact-mdd: Efficiently filtering (s)mdd constraints with
  reversible sparse bit-sets\BBCQ\
\newblock In {\Bem Proceedings of the Twenty-Seventh International Joint
  Conference on Artificial Intelligence}, \BPGS\ 1383--1389.

\bibitem[\protect\BCAY{Vion\ \BBA\ Piechowiak}{Vion\ \BBA\
  Piechowiak}{2018}]{global3}
Vion, J.\BBACOMMA\  \BBA\ Piechowiak, S. \BBOP2018\BBCP.
\newblock \BBOQ From mdd to bdd and arc consistency\BBCQ\
\newblock {\Bem Constraints}, {\Bem 23}.

\end{thebibliography}

\phantom{\cite{tsptwbench,MDDForSP,ddo,BnB,BnB2,hadzic2008,bergman2011,andersen2007,hoda2010,jobsequencing,canonicalarcs,tsplib,Perez_Regin_2018,Rudich2022PeelAndBoundGS,tsptwbench,vianney,baldacci2012new,gillard2022discrete,bergman2014optimization}}
\end{document}